\title{Giant components in biased graph processes}
\date{}
\author{{Gideon Amir\thanks{Department of Mathematics,
University of Toronto, Toronto, ON M5S-2E4, Canada. Email address: {\tt gidi@math.utoronto.ca}.}}
\quad {Ori Gurel-Gurevich\thanks{Theory Group of Microsoft Research, One Microsoft Way, Redmond,
WA 98052-6399, USA. Email address: {\tt origurel@microsoft.com}.}} \quad {Eyal Lubetzky\thanks{Theory Group of Microsoft Research, One Microsoft Way, Redmond,
WA 98052-6399, USA. Email address: {\tt eyal@microsoft.com}.}}
\quad {Amit Singer\thanks{Department of Mathematics and PACM, Princeton University, Princeton, NJ 08544-1000, USA. Email address: {\tt amits@math.princeton.edu}.}}}

\documentclass [11pt]{article}
\usepackage[]{graphicx,epsf,epsfig,amsmath,amssymb,amsfonts,latexsym,amsthm}

\newtheorem{theorem}{Theorem}[section]
\newtheorem{lemma}[theorem]{Lemma}
\newtheorem{claim}[theorem]{Claim}
\newtheorem*{definition}{Definition}

\newtheorem{corollary}[theorem]{Corollary}
\newtheorem{proposition}[theorem]{Proposition}

\renewcommand{\epsilon}{\varepsilon}

\newcommand{\Gorg}[1][K]{{\mathcal{G}_#1}}
\newcommand{\Gapp}[1][K]{{\widetilde{\mathcal{G}_#1}}}
\newtheoremstyle{upright}%
        {8pt plus2pt minus4pt}%
        {8pt plus2pt minus4pt}%
        {\upshape}%
        {}%
        {\bfseries\scshape}%
        {:}%
        {1em}%
        {}%
\theoremstyle{upright}
\newtheorem{remark}[theorem]{Remark}
\newcommand{\ds}{\displaystyle}
\newcommand{\ignore}[1]{}

\begin{document}
\maketitle

\begin{abstract}
A random graph process, $\Gorg[1](n)$, is a sequence of graphs on
$n$ vertices which begins with the edgeless graph, and where at each
step a single edge is added according to a uniform distribution on
the missing edges. It is well known that in such a process a giant
component (of linear size) typically emerges after
$(1+o(1))\frac{n}{2}$ edges (a phenomenon known as ``the double
jump''), i.e., at time $t=1$ when using a timescale of $n/2$ edges
in each step.

We consider a generalization of this process, $\Gorg(n)$, proposed
by Itai Benjamini in order to model the spreading of an epidemic. This generalized
process gives a weight of size $1$ to missing edges between pairs of
isolated vertices, and a weight of size $K \in [0,\infty)$
otherwise. This corresponds to a case where links are added between
$n$ initially isolated settlements, where the probability of a new
link in each step is biased according to whether or not its two
endpoint settlements are still isolated.

Combining methods of \cite{SpencerWormald} with analytical
techniques, we describe the typical emerging time of a giant
component in this process, $t_c(K)$, as the singularity point of a
solution to a set of differential equations. We proceed to analyze
these differential equations and obtain properties of $\Gorg$, and
in particular, we show that $t_c(K)$ strictly decreases from
$\frac{3}{2}$ to $0$ as $K$ increases from $0$ to $\infty$, and that
$t_c(K) = \frac{4}{\sqrt{3K}}\left(1 + o(1)\right)$, where the $o(1)$-term
tends to $0$ as $K\to\infty$. Numerical
approximations of the differential equations agree both with
computer simulations of the process $\Gorg(n)$ and with the
analytical results.
\end{abstract}

\section{Introduction}\label{sec::intro}
\subsection{The Achlioptas problem and the biased process}
The random graph process on $n$ vertices, $\Gorg[1]=\Gorg[1](n)$,
introduced by Erd\H{o}s and R\'enyi, is a sequence of
$\binom{n}{2}+1$ graphs,
$(\Gorg[1]^0,\ldots,\Gorg[1]^{\binom{n}{2}})$, where the
$\Gorg[1]^0$ is the edgeless graph on $n$ vertices, and $\Gorg[1]^m$
is obtained by adding a random edge to $\Gorg[1]^{m-1}$, chosen
uniformly over all missing edges.

A classical result of Erd\H{o}s and R\'enyi (\cite{ErdosRenyi})
states that if $T= C \frac{n}{2} $ and $C<1$, then typically every
connected component of $\Gorg[1]^T$ is of size $O(\log n)$, and if
$C>1$ then typically there is a single {\em giant component} of size
$\Theta(n)$ and every other component is of size $O(\log n)$. Thus,
there is a phase transition (the ``double jump'') after
$(1+o(1))n/2$ edges. For further discussion of this phenomenon, see,
e.g., \cite{RandomGraphs}.

A well known problem, introduced by Achlioptas, discusses a scenario
where two randomly chosen edges are presented at each step, out of
which a single edge is chosen by some algorithm $\mathcal{A}$. The
goal of the algorithm is to postpone the emerging time of the giant
component as much as possible. This was first examined by Bohman and
Frieze in \cite{BohmanFrieze1}. For additional results both on this
problem (as an off-line and as an online problem), as well as on the
converse problem of creating a giant component ahead of time, see
\cite{BohmanFrieze2},\cite{BohmanKim},\cite{BohmanKravitz},\cite{FlaxmanGamarnikSorkin}.
In \cite{SpencerWormald}, the authors describe a generic approach to
analyzing the performance of algorithms for the mentioned Achlioptas
problem. After applying Wormald's differential equation method for
graph processes \cite{WormaldDiffEq}, the emerging time of the giant
component is expressed as a singularity point to a differential
equation. Using this method, the authors are able to provide bounds
for the performance of several algorithms.

In this paper, we study a natural generalization of the Erd\H{o}s-R\'enyi
random graph process, proposed by Itai Benjamini
in order to model the spreading of an epidemic. This process of Benjamini lets ``infected'' clusters
have either a larger or a smaller probability of increasing their size, depending
on the value of an external continuous parameter $K$. This is achieved by embedding basic degree information (namely, whether or not a site is currently isolated) into the probability distribution over the missing edges. Note that the given model is the most natural of its kind with respect to the dependency of this probability distribution on the degrees of the vertices. Results on other, more complicated, models may be obtained via methods similar to the ones presented here.

We combine
methods of \cite{SpencerWormald} together with
analytical methods, with the same motivation of determining the
critical time in which a giant component emerges in this process. 
As we mentioned, the generalized
process we consider is a parameterized version, dependent on some
$K\in[0,\infty)$, which modifies the probability of each edge
according to whether or not its endpoints are already connected. In
the original Achlioptas problem, different algorithms can postpone
the phase transition or create it ahead of time, where the biased
choice at each step is between precisely two randomly chosen edges.
The process we study considers all missing edges when making its
biased choice, and the phase transition is presented as a function
of the continuous parameter $K$.

By applying the powerful differential equation method, we were able
to derive properties of our model directly from the system of
coupled non-linear ordinary differential equations (ODEs). While
some of these properties can be proved by relatively simple
combinatorial arguments, calculating the precise asymptotic behavior
of $t_c(K)$ for $K \gg 1$ via combinatorial arguments seems
challenging. Indeed, while at first glance this result appears as though
it can be obtained by probabilistic arguments (e.g., using monotone coupling to variants of the model, tracking the structure of the connected components throughout certain time intervals, etc.), we are
not aware of any such derivation at the present time.
In particular, the only way to obtain the $4/\sqrt{3K}$ behavior of the blowup time is through the careful asymptotic analysis of the derived coupled differential equations
(see, e.g., \cite{Bender} for more on asymptotic analysis of ODEs).

The generalized process can be efficiently implemented. In order to
efficiently randomize the next edge, one needs to maintain the sets
of isolated and non-isolated vertices, along with the set of edges
already chosen. Our implementation runs in time $O(n \log n)$ and
requires $O(n)$ memory. The computer simulations show an excellent
agreement with the numerical solutions of the ODEs on the one hand,
and with the analytical results concerning $t_c(K)$ on the other
hand.

\subsection{Notations and main results}
The {\em biased graph process} on $n$ vertices, $\Gorg$, is
the following generalization of the random graph process: as before,
the initial graph $\Gorg^0$ is the edgeless graph on $n$ vertices,
and $\Gorg^m$ is obtained by adding a single edge to $\Gorg^{m-1}$.
The newly added edge is selected according to the following
distribution on the missing edges: each edge between two isolated
vertices is assigned a weight of $1$, and the weight $K$ is assigned
to all the remaining edges. Once there are less than $2$ isolated
vertices, the distribution on the missing edges is uniform. We
extend the definition of $\Gorg^m$ to $m
> \binom{n}{2}$ by setting $\Gorg^m=K_n$ for every such $m$, where
$K_n$ denotes the complete graph on $n$ vertices. Furthermore, we
use the notation $\Gorg^T|_H$, where $H$ is a graph on $n$ vertices,
to denote the biased graph process after $T$ steps, starting from
the initial graph $H$ instead of the edgeless graph.

Let $\Gorg(t)$ denote the biased process after scaling its time line
by a factor of $n/2$:
$$\Gorg(t) = \Gorg^{\lfloor t n / 2 \rfloor}~.$$
Since a choice of $K=1$ is equivalent to the Erd\H{o}s-R\'enyi
random graph process $\Gorg[1]$, the appearing time of the giant
component in $\Gorg[1]$ is typically at $t=1$. We study the effect
that modifying $K$ has on this critical point, $t_c(K)$, keeping in
mind that, intuitively, decreasing the value of $K$ should postpone the
emerging point of a giant component and vice versa.

Throughout the paper, we say that a random graph on $n$ vertices
satisfies some property {\em with high probability}, or {\em almost
surely}, or that {\em almost every} graph process on $n$ vertices
satisfies a property, if the corresponding event has a probability
which tends to $1$ as $n$ tends to infinity.

Let $G=(V,E)$ be a graph on $|V|=n$ vertices. Let
$\mathcal{C}=\mathcal{C}(G)$ denote the set of connected components
of $G$, and let $\mathcal{C}_v=\mathcal{C}_v(G)$ denote the
connected component of a given vertex $v \in V$. We denote by
$\mathcal{C}_0$ the set of isolated vertices in $G$, and by $I(G)$
the fraction of isolated vertices:
$$I(G) = \frac{|\mathcal{C}_0|}{n}~.$$
Notice that, if the graph after $T$ steps (that is, once $T$ edges have been added) contains $|\mathcal{C}_0|$ isolated vertices, then the
probability that the next edge is a specific one between two
isolated vertices is $$\bigg[\binom{|\mathcal{C}_0|}{2} + \left(\binom{n}{2}-\binom{|\mathcal{C}_0|}{2}-T\right)K\bigg]^{-1}~,$$ and the probability that it is a specific new edge between two vertices that are not both isolated is $K$ times that value.

The {\em susceptibility} of $G$, $S(G)$, is defined to be the
expected size of a connected component of a uniformly chosen vertex
$v \in V$:
$$S(G) = \frac{1}{n}\sum_{v \in V}|\mathcal{C}_v| =
\frac{1}{n}\sum_{C \in \mathcal{C}(G)} |C|^2~.$$ The susceptibility
of $G$ is closely coupled with the existence of a giant component in
$G$. Indeed, the existence of a giant component of size $\alpha n$
for some $\alpha > 0$ implies that $S(G) \geq \alpha^2 n$, and
conversely, if $S(G) = \Omega(n)$ then at least one connected
component is of linear size. Following the ideas of
\cite{SpencerWormald}, we characterize the behavior of the
susceptibility along $\Gorg$, and in the process obtain the required
results on $t_c(K)$. In order to do so, we need to examine the
typical behavior of the number of isolated vertices along $\Gorg$.

The fraction of isolated vertices in $\Gorg$, $I(\Gorg(t))$ has a
value of $1$ at $t=0$, and decreases to $0$ over time. The following
theorem summarizes the behavior of $I(\Gorg(t))$:
\begin{theorem}\label{I-thm}
For every $K>0$ and every $C>0$, almost every biased process on $n$
vertices $\Gorg$ satisfies $|I(\Gorg(t))-y(t)|=O(n^{-1/4})$ for
every $0 \leq t \leq C$, where $y$ is the solution to the
differential equation:
\begin{equation}\label{y-eq}
  \left\{\begin{array}{lcl}
y' &=& \displaystyle{\frac{(1-y)K}{y^2+(1-y^2)K}-1} \\
y(0)&=&1 \\
\end{array} \right. ~.
\end{equation}
In the special case $K=0$, the above holds for $C=1$.
\end{theorem}

The following proposition analyzes the differential equation whose
solution will prove to be a good approximation of the susceptibility
along the biased process:

\begin{proposition}\label{z-thm}
Let $z=z(t)$ denote the solution for the following differential
equation:
\begin{equation}\label{z-eq}
  \left\{\begin{array}{lcl}
z' &=& \ds{\frac{K}{y^2+(1-y^2)K} \left(z^2 - 1\right) + 1}\\
z(0)&=&1 \\
\end{array} \right. ~,
\end{equation}
where $y$ is the solution to the differential equation \eqref{y-eq}.
For every $K>0$ there exists a singularity point $t_c = t_c(K) > 0$,
such that $z(t)$ is continuous on $[0,t_c)$, and $\lim_{t\rightarrow
t_c^-}z(t)=\infty$. Furthermore, there exists some constant $M > 0$,
independent of $K$, such that $t_c(K) \in (0,M)$ for every $K>0$.
\end{proposition}

The next theorem implies that the singularity point of $z(t)$,
$t_c=t_c(K)$, is the typical time at which a giant component emerges
in a biased process:

\begin{theorem}\label{S-thm}
For every $\epsilon > 0$ and every $K>0$, almost every biased
process on $n$ vertices $\Gorg$ satisfies the following:
\begin{enumerate}
\item \label{S-thm-sub} Subcritical phase: $|S(\Gorg(t)) - z(t)| = o(1)$ for every $t\in [0,t_c-\epsilon]$,
where $z(t)$ is the solution for the differential equation
\eqref{z-eq}, and $t_c$ is its singularity point as defined in
Proposition \ref{z-thm}. Furthermore, for every $t \in
[0,t_c-\epsilon]$, the largest component of $\Gorg(t)$ is of size
$O(\log n)$.
\item \label{S-thm-super} Supercritical phase: $S(\Gorg(t_c + \epsilon)) = \Omega(n)$.
\end{enumerate}
Altogether, the appearance of a giant component in $\Gorg$ is almost
surely at time $t_c$. In the special case $K=0$, the above holds
when replacing $z(t)$ with the function $\hat{z}(t)=
\left\{\begin{array}{cl}
1+t & \text{if } t\leq 1 \\
\frac{1}{3/2-t} & \text{if } t \geq 1 \\
\end{array} \right.$.
\end{theorem}

Finally, the behavior of $t_c(K)$, the typical point of the phase
transition in $\Gorg$, is characterized by the following theorem:

\begin{theorem}\label{critical-thm}
Let $t_c(K)$ denote the singularity point of the solution to the ODE
\eqref{z-eq}, as defined in Proposition \ref{z-thm}. Then $t_c(K)$
is continuous and strictly monotone decreasing as a function of $K$
and satisfies:
$$\left\{\begin{array}{lcl}
t_c(0)&=&\frac{3}{2}\\
t_c(K)&=&\frac{4}{\sqrt{3K}}\left(1 + o(1)\right)\\
\end{array} \right.~,$$
where the $o(1)$-term tends to $0$ as $K\to\infty$.
\end{theorem}

The rest of the paper is organized as follows: in Section
\ref{sec::distributions} we study general properties of the
distribution of $\Gorg$. Namely, we study the relation between
$\Gorg$ and $\Gorg[1]$, and describe an approximated process,
$\Gapp$, which is easier to analyze.

In Section \ref{sec::isolation} we prove Theorem \ref{I-thm} and
analyze the solution to equation \eqref{y-eq}, which characterizes
the behavior of the isolation ratio throughout the biased process.
Proposition \ref{z-thm} and Theorem \ref{S-thm} are both proved in
Section \ref{sec::susceptibility}, where computer simulations are
also included.

In Section \ref{sec::critical} we analyze the asymptotic behavior of
$t_c(K)$ for large values of $K$, and prove Theorem
\ref{critical-thm}. The final section \ref{sec::conclusion} is
devoted to open problems and concluding remarks.

\section{Dominating and approximate distributions for $\Gorg$}
\label{sec::distributions}
\subsection{The relation between $\Gorg$ and $\Gorg[1]$}

A main tool in understanding and analyzing the biased process
$\Gorg$ is the stochastic domination relation between $\Gorg$ and
time-stretched versions of the Erd\H{o}s-R\'enyi process $\Gorg[1]$.
In order to formalize and prove this relation, we consider a wider
family of graph processes, defined as follows:

\begin{definition} Let $M \in \mathbb{N}$.
An \texttt{$M$-bounded weighted graph process} on $n$ vertices,
$\mathcal{H}=\mathcal{H}(n)$, is an infinite sequence of graphs on
$n$ vertices, $(\mathcal{H}^0,\mathcal{H}^1,\ldots)$, where
$\mathcal{H}^0$ is some fixed initial graph, and $\mathcal{H}^t$ is
generated from $\mathcal{H}^{t-1}$ by adding one edge at random, as
follows: the probability of adding the edge $e$ to
$\mathcal{H}^{t-1}$ is proportional to some weight function
$W_t(e)$, satisfying:
$$\max_{e\notin \mathcal{H}^{t-1}} W_t(e) \leq M \min_{e\notin H^{t-1}}
W_t(e)~.
$$If for some $\nu\geq 0$ $\mathcal{H}^\nu = K_n$, we define
$\mathcal{H}^t = \mathcal{H}^\nu = K_n$ for every $t>\nu$.
\end{definition}
Clearly, the biased process $\Gorg$ is an $M$-bounded weighted graph
process which starts from the edgeless graph, where $M=\lceil
\max\{K,\frac{1}{K}\}\rceil$. We are interested in the relation
between the probability that $\Gorg[1]$ satisfies some monotone
graph property $\mathcal{A}$ (a set of graphs closed under
isomorphism and under the addition of edges) and the corresponding
probability of $\Gorg$. The following theorem formalizes the
stochastic domination of the original graph process $\Gorg[1]$ on
$M$-bounded weighted graph processes:
\begin{theorem}\label{thm-M-bounded}
Let $\mathcal{H}$ denote an $M$-bounded weighted graph process on
$n$ vertices, and let $\mathcal{A}$ denote a monotone increasing
property of graphs on $n$ vertices. The following statements hold
for any $t \in \mathbb{N}$:
\begin{eqnarray}
\label{dom-M-by-1} \Pr[\mathcal{H}^t \in \mathcal{A}]
&\leq& \Pr[\Gorg[1]^{M t}|_{\mathcal{H}^0} \in \mathcal{A} ] ~,\\
\label{dom-1-by-M} \Pr[\Gorg[1]^{t}|_{\mathcal{H}^0} \in
\mathcal{A}]&\leq& \Pr[\mathcal{H}^{M t}\in \mathcal{A}]~.
\end{eqnarray}
\end{theorem}

We need the following lemma, which was first proved in
\cite{Strassen} in a slightly different setting. For the sake of
completeness, we include a short proof of the lemma using the
Max-Flow-Min-Cut Theorem (a relation which was first observed in
\cite{Preston}):
\begin{lemma}\label{flow-lemma}
Let $U,V$ be two finite sets, and let $R\subset U \times V$ denote a
relation on $U,V$. Let $\mu$ and $\nu$ denote probability measures
on $U$ and $V$ respectively, such that the following inequality
holds for every $A \subset U$:
\begin{equation}\label{hall-condition}
\mu(A) \leq \nu(\{y\in V : x R y \mbox{ for some }x \in A\})~.
\end{equation}
Then there exists a coupling $\varphi$ of $\mu,\nu$ whose support is
contained in $R$. That is, there is a joint distribution $\varphi$
on $U,V$ satisfying the following two properties:
\begin{enumerate}
  \item \label{lemma-marginal-property}
  The marginal distributions of $\varphi$ on $A$ and $B$
are $\mu$ and $\nu$ respectively.
  \item \label{lemma-inclusion-property} For every $(x,y)\in U \times V$ such that
  $\neg(x R y)$, $\varphi(x,y)=0$.
\end{enumerate}
Furthermore, the distribution $\varphi$ can be found in time
polynomial in $|U|+|V|$.
\end{lemma}
\begin{figure}
\centering
\includegraphics{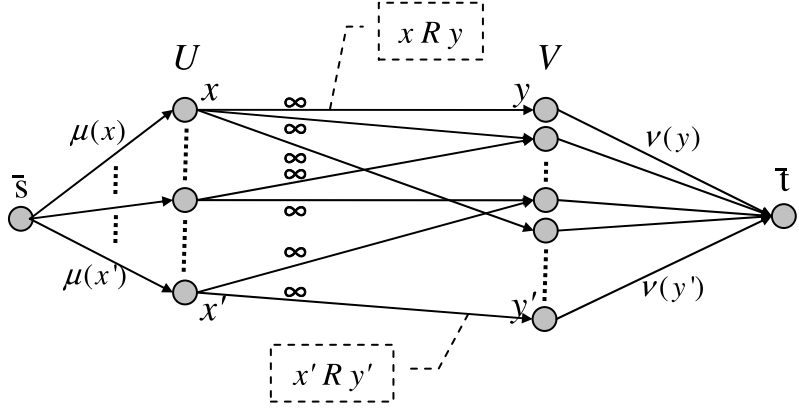}
\caption{The network graph $G$ in Lemma \ref{flow-lemma}.}
\label{fig::network-flow}
\end{figure}
\begin{proof}
Consider a directed weighted graph $G$ on the vertex set $U\cup
V\cup\{\bar{s},\bar{t}\}$, with the following set of edges and
capacities (see Figure \ref{fig::network-flow}):
\begin{enumerate}
\item For every $x \in U$ and $y \in V$ such that $x R y$, place
the edge $(x,y)$ with a capacity of $\infty$.
\item For every $x \in U$, place the edge $(\bar{s},x)$ with
a capacity of $\mu(x)$.
\item For every $y \in V$, place the edge $(y,\bar{t})$ with a
capacity of $\nu(y)$.
\end{enumerate}
Note that a flow of $1$ from $\bar{s}$ to $\bar{t}$ in the network
defined above provides the desired distribution $\varphi$, by
defining $\varphi(x,y)=f(x,y)$ where $f$ is such a flow. The flow
value of $1$ implies the first property required from $\varphi$. The
second property follows from the fact that, if $\neg(x R y)$ for
some $x \in U$ and $y \in V$, then the edge $(x,y)$ is not in $G$.
Altogether, if indeed $G$ has a flow of $1$, then a joint
distribution $\varphi$ satisfying the above properties can be
computed in time polynomial in $|U|+|V|$.

The fact $G$ has a flow of $1$ follows essentially from the proof of
Hall's Theorem using the Max-Flow-Min-Cut Theorem (see, e.g.,
\cite{MatchingTheory}), where inequality \eqref{hall-condition}
replaces Hall's criteria for a maximal matching. Let $f(G)$ denote
the maximal flow from $\bar{s}$ to $\bar{t}$ in $G$. By the
Max-Flow-Min-Cut Theorem, $f$ corresponds to a minimal $(S,T)$ cut,
where $\bar{s}\in S$ and $\bar{t}\in T$. Hence, if we define $A =
S\cap U$ and $B= S \cap V$, the choice of the capacity $\infty$ for
the edges between $U,V$ implies that $N(A) \cap V \subset B$, where
$N(A)$ is the set of neighbors of $A$ in the graph $G$. Hence:
$$ f = \mu(U \setminus A) + \nu(B) \geq \mu(U \setminus A) + \nu(N(A) \cap V) \geq
\mu(U \setminus A) + \mu(A) = 1~,$$ where the last inequality is by
the assumption $\eqref{hall-condition}$.
\end{proof}

\begin{proof}[Proof of Theorem \ref{thm-M-bounded}] Let
$\mathcal{P},\mathcal{Q}$ denote two graph processes starting from
the same initial graph. We wish to prove that, under certain
conditions, the following inequality holds for every $t \in
\mathbb{N}$:
\begin{equation}\label{eq-p-q-domination} \Pr[
\mathcal{P}^t \in \mathcal{A}] \leq \Pr[\mathcal{Q}^{M t} \in
\mathcal{A}] ~.
\end{equation} The proof will follow from a
coupling of the two processes $\mathcal{P}$ and $\mathcal{Q}$, such
that for every instance of the coupling and every $t\in\mathbb{N}$,
$\mathcal{P}^t$ is a subset of $\mathcal{Q}^{t M}$. More precisely,
we define a joint distribution on the processes $
(\mathcal{P},\mathcal{Q})$ in the following manner: at each step $t
\in \mathbb{N}$, we generate $\mathcal{P}^t$ according to its usual
law, then generate the pairs $\mathcal{Q}^{(t-1) M +
1},\ldots,\mathcal{Q}^{t M}$ according to some distribution, such
that $E(\mathcal{P}^t) \subset E(\mathcal{Q}^{t M})$ and the
marginal distribution of $\mathcal{Q}^{(t-1) M +
1},\ldots,\mathcal{Q}^{t M}$ is the correct one. Clearly, such a
construction immediately proves \eqref{eq-p-q-domination}.

For every $t \in \mathbb{N}$ and $i\in\{1,\ldots,M\}$, let $X_t^i$
and $Y_t$ denote the sets of missing edges of
$\mathcal{Q}^{(t-1)M+i-1}$ and $\mathcal{P}^{t-1}$ respectively. Let
$e_t \in Y_t$ denote the $t$-th edge added to $\mathcal{P}$, and let
$f_t^i \in X_t^i$ denote the $((t-1)M+i)$-th edge added to
$\mathcal{Q}$.

Recalling that, by definition, $\mathcal{P}^0 = \mathcal{Q}^0$, let
$t \geq 1$ and assume that we have constructed the above coupling up
to the point $(\mathcal{P}^{t-1},\mathcal{Q}^{(t-1)M})$. In
particular, $E(\mathcal{P}^{t-1})\subset E(\mathcal{Q}^{(t-1)M})$,
and hence $X_t^1 \subset Y_t$.

Note that in the special case $M \geq |X_t^1|$, in which
$\mathcal{Q}^{t M}$ is by definition the complete graph, we can
clearly generate
$\mathcal{Q}^{(t-1)M+1},\ldots,\mathcal{Q}^{(t-1)M+|X|}$ according
to the usual law of $\mathcal{Q}$ and maintain the properties
required from the coupling. Assume therefore that $M < |X_t^1|$.

We now wish to extend the coupling to $(\mathcal{P}^{t
M},\mathcal{Q}^t)$ using the joint distribution $\varphi$ of Lemma
\ref{flow-lemma}. To do so, take:
\begin{equation}\label{values-for-flow-lemma}\left\{\begin{array}{lcl}
U&=&Y_t\\
V&=&(X_t^1)_M\\
\mu(e)&=& \Pr[e_t=e] ~\mbox{ for every }e \in Y_t\\
\nu(F)&=&\Pr[(f_t^1,\ldots,f_t^M)=F]~\mbox{ for every }F \in
(X_t^1)_M\\
R &=& \{(e,F) \in Y_t \times (X_t^1)_M ~:~e \in (Y_t \setminus
X_t^1) \mbox{ or } e \in F\}
\end{array} \right. ~,
\end{equation}
where $(X_t^1)_M$ denotes the set of ordered subsets of $M$ edges of
$X_t^1$.

If inequality \eqref{hall-condition} holds for every $A \subset U$,
then $\varphi$ associates every edge $e \in Y_t$ with a distribution
on $\mathcal{Q}^{(t-1)M+1},\ldots,\mathcal{Q}^{t M}$, which ensures
that $e \in \mathcal{Q}^{t M}$, by the second property of the lemma.
The choice of $\nu$, along with the first property of the lemma,
ensures that $\mathcal{Q}^{(t-1)M+1},\ldots,\mathcal{Q}^{t M}$ will
have the correct marginal distributions.

Clearly, \eqref{hall-condition} holds for every $A \subset U$ such
that $A \cap (Y_t \setminus X_t^1) \neq \emptyset$, since by the
definition of $R$, the right hand side of $\eqref{hall-condition}$
equals $1$ in this case. It remains to show that inequality
\eqref{hall-condition} indeed holds for every $A \subset X_t^1$ when
$\mathcal{P}$ and $\mathcal{Q}$ play the roles of the processes
$\Gorg[1]|_{\mathcal{H}^0}$ and $\mathcal{H}$. This is proved in the
following two claims:

\begin{claim}\label{claim-hall-condition-1-holds}
Let $\mathcal{P}=\mathcal{H}$ and
$\mathcal{Q}=\Gorg[1]|_{\mathcal{H}^0}$, and define $U,V,\mu,\nu,R$
as in \eqref{values-for-flow-lemma}. Then inequality
\eqref{hall-condition} holds for every $A \subset U$.
\end{claim}
\begin{proof} To simplify the notations, let $X=X_t^1$ and
let $Y=Y_t$. By the above choices for $\mathcal{P}$ and
$\mathcal{Q}$ and by the definition of the $M$-bounded weighted
process, we have:
$$\nu(F) = \Pr[(f_t^1,\ldots,f_t^M)=F] = \frac{1}{|(X)_M|}~\mbox{ for every }F\in(X)_M~,$$
\begin{equation}\label{m-bounded-edge-probability}
\mu(e)=\Pr[e_t = e] = \frac{W_t(e)}{\sum_{e' \in Y}W_t(e')} ~\mbox{
for every }e\in Y~,
\end{equation}
where, without loss of generality:
$$\min_{e \in Y}W_t(e)=1 ~,~ \max_{e \in Y}W_t(e)\leq M~.$$
Substituting the values of $\mu,\nu$, inequality
\eqref{hall-condition} takes the form: for every $A \subset X$,
\begin{equation}\label{hall-condition-1}
\Pr[e_t \in A] \leq 1 - \frac{\binom{|X|-|A|}{M}}{\binom{|X|}{M}}~.
\end{equation}
Take $A \subset X$, and assume that $A \neq \emptyset,X$ (otherwise
\eqref{hall-condition-1} trivially holds). By
\eqref{m-bounded-edge-probability}, we have:
$$\Pr[e_t \in A] = \sum_{e \in A}\Pr[e_t = e] = \frac{\sum_{e\in
A}W_t(e)}{\sum_{e\in Y}W_t(e)} \leq \frac{|A|M}{|Y|+(M-1)|A|}~,$$
where the last inequality is by the fact that $1 \leq W_t(e) \leq M$
for every $e\in Y$, along with the inequality
$\frac{a+c}{b+c}>\frac{a}{b}$ for every $0<a<b$ and $c>0$. Recalling
that $X \subset Y$, we get:
\begin{equation}\label{hall-lhs-bound} \Pr[e_t \in
A] \leq \frac{|A|M}{|X|+(M-1)|A|}~.
\end{equation} Combining \eqref{hall-lhs-bound} with the fact that:
\begin{equation}\label{hall-rhs-bound}
1-\frac{\binom{|X|-|A|}{M}}{\binom{|X|}{M}} = 1 - \prod_{j=0}^{M-1}
\frac{|X|-|A|-j}{|X|-j} \geq 1 - \left(1-\frac{|A|}{|X|}\right)^M~,
\end{equation}
we obtain the following sufficient condition for
\eqref{hall-condition-1}:
\begin{equation}\label{suf-condition-for-hall-1}
\left(1-\frac{|A|}{|X|}\right)^M \leq \frac{|X|-|A|}{|X|-|A|+M|A|}~.
\end{equation}
Set $\beta = 1-\frac{|A|}{|X|}$, and recall that $0<\beta<1$.
Inequality \eqref{suf-condition-for-hall-1} takes the following
form:
$$ \beta^M \leq \frac{\beta}{\beta+(1-\beta)M}~,$$
or equivalently:
\begin{equation}\label{suf-condition-for-hall-2}
\beta^{1-M}-(1-M)\beta-M \geq 0~.
\end{equation} Defining $g(x) = x^{1-M}-(1-M)x-M$, it is easy to
verify that $g'(x)\leq0$ for $0<x<1$ provided that $M\geq 1$, and
that $g(1)=0$. Hence, inequality \eqref{suf-condition-for-hall-2}
indeed holds, and inequality \eqref{hall-condition-1} follows, as
required.
\end{proof}

\begin{claim}\label{claim-hall-condition-2-holds}
Let $\mathcal{P}=\Gorg[1]|_{\mathcal{H}^0}$ and
$\mathcal{Q}=\mathcal{H}$, and define $U,V,\mu,\nu,R$ as in
\eqref{values-for-flow-lemma}. Then inequality
\eqref{hall-condition} holds for every $A \subset U$.
\end{claim}
\begin{proof} Following the notation of the previous claim,
define $X=X_t^1$ and $Y=Y_t$. As $\mu$ is the uniform distribution
on $Y$, inequality \eqref{hall-condition} takes the following form:
for every $A \subset X$,
\begin{equation}\label{hall-condition-2}
\frac{|A|}{|Y|} \leq 1-\Pr[f_t^1,\ldots,f_t^M \notin A]~.
\end{equation}
Let $A \subset X$, and we may again assume $A \neq \emptyset,X$
otherwise \eqref{hall-condition-2} trivially holds. For
$i=1,\ldots,M$, let $W_t^i$ denote the weight function
$W_{(t-1)M+i}$ on the set of edges $X_t^i$. According to this
notation, the right hand side of \eqref{hall-condition-2} satisfies:
$$1 -
\Pr[f_t^1,\ldots,f_t^m \notin A] = 1 - \prod_{i=1}^M \Pr[f_t^i
\notin A ~|~ f_t^1,\ldots,f_t^{i-1} \notin A] = $$ $$ = 1 -
\prod_{i=1}^M \left(1-\ds{\frac{\sum_{e \in A} W_t^i(e) }{\sum_{e
 \in X_t^i} W_t^i(e) }}\right)~.$$
By the same argument used in Claim
\ref{claim-hall-condition-1-holds}, we reduce the expression
$\ds{\frac{\sum_{e \in A} W_t^i(e) }{\sum_{e
 \in X_t^i} W_t^i(e) }}$ by assigning the value $1$ to the weights of
 $A$ and $M$ to the rest. Hence,
$$1 - \Pr[f_t^1,\ldots,f_t^m \notin A] \geq
 1 - \prod_{i=1}^M
\left(1-\frac{|A|}{|X_t^i|M + |A|(1-M)}\right)\geq 1 -
\left(\frac{|X|M-|A|M}{|X|M-|A|M+|A|}\right)^M~.$$ Recalling that $X
\subset Y$ and setting $\beta =1-\frac{|A|}{|X|}$ ($0<\beta<1$), the
claim then follows from the next inequality:
\begin{equation}
  1 - \beta \leq 1 - \left(\frac{\beta}{\beta M +
  (1-\beta)}\right)^M~,
\end{equation}
which holds for every $\beta\geq 0$ whenever $M\geq1$. \end{proof}
Combining Claim \ref{claim-hall-condition-1-holds} with the
arguments preceding it and Lemma \ref{flow-lemma}, yields that it is
possible to extend the coupling to $(\mathcal{P}^{t
M},\mathcal{Q}^t)$ when $\mathcal{P}=\mathcal{H}$ and
$\mathcal{Q}=\Gorg[1]|_{\mathcal{H}^0}$, thereby completing the
induction argument. Therefore, \eqref{dom-M-by-1} holds for every $t
\in \mathbb{N}$. Similarly, combining Claim
\ref{claim-hall-condition-2-holds} with the above arguments yields
that \eqref{dom-1-by-M} holds for every $t\in\mathbb{N}$. This
completes the proof of the theorem.
\end{proof}

\begin{remark}
Theorem \ref{thm-M-bounded} required that one of the processes,
$\mathcal{H}$, chooses each edge according to an $M$-bounded
distribution, whereas the other chooses each edge according to a
uniform distribution. It is not difficult to construct an example
showing that this requirement cannot be replaced by the condition,
that the maximal ratio between the weights of the two processes at
each step is at most $M$.
\end{remark}

Applying Theorem \ref{thm-M-bounded} on the biased graph process
$\Gorg$ gives the following immediate corollary:
\begin{corollary}\label{cor-stochastic-domination} Let $\mathcal{A}$
denote a monotone increasing property of graphs on $n$ vertices, and
let $t\geq 0$. Then for every $K>0$, the following two statements
hold:
\begin{enumerate}
\item If $\Gorg(t)$ almost surely satisfies $\mathcal{A}$, then
$\Gorg[1](\lceil\max\{K,\frac{1}{K}\}\rceil t)$ almost surely
satisfies $\mathcal{A}$.

\item If $\Gorg[1](t)$ almost surely satisfies $\mathcal{A}$, then
$\Gorg(\lceil\max\{K,\frac{1}{K}\}\rceil t)$ almost surely satisfies
$\mathcal{A}$.
\end{enumerate}
\end{corollary}

\subsection{The approximate biased graph process}
In order to simplify the proofs of Theorems \ref{I-thm} and
\ref{S-thm}, we consider a variant of the biased process on $n$
vertices, $\Gapp$, which we dub an {\em approximate}
biased process. At each step, a random ordered pair of vertices
$(u,v)\in V^2$ (where $V$ is the set of vertices) is chosen out of
the $n^2$ possible pairs, according to the following distribution: a
pair of isolated vertices has a weight of $1$, whereas all other
pairs have weights of $K$. If the chosen pair corresponds to a
(self) loop or to an edge which already exists in $\Gapp$, no edge
is added in this step.

The following claim implies that it is sufficient to prove Theorems
\ref{I-thm} and \ref{S-thm} for the approximate model:
\begin{claim}\label{app-process-clm}
Let $C > 0$, and let $\{\mathcal{A}_t : 0 \leq t \leq C\}$ denote a
family of properties of graphs on $n$ vertices. If for almost every
approximate biased process $\Gapp$, $\Gapp(t)$ satisfies
$\mathcal{A}_t$ for every $0\leq t \leq T$, then for almost every
biased process $\Gorg$, $\Gorg(t)$ satisfies $\mathcal{A}_t$ for
every $0 \leq t \leq T$ as well.\end{claim}
\begin{proof}
Let $B_j$ denote the event that the ordered pair, chosen in the
$j$-th step, was not added to $\Gapp$, either being a loop or
already belonging to $\Gapp$. Fix $C > 0$, and let $T=C
\frac{n}{2}$. The probability of the event $B_j$ ($1\leq j \leq T$)
satisfies:
$$\Pr[B_j] \leq \frac{1}{n} + \frac{2(j-1)}{n^2} \leq \frac{C+1}{n}~,$$
and, as these events are independent, we apply the well known bound
$1-x\geq \mathrm{e}^{-x/(1-x)}$ for $0 \leq x<1$, and obtain:
$$\Pr[\wedge_{j=1}^T \overline{B_i}] \geq
\left(1-\frac{C+1}{n}\right)^T \geq
\exp\left(-\frac{\frac{C+1}{n}T}{1-o(1)}\right) \geq
\exp(-C(C+1))~.$$ Notice that if we condition on the event
$\wedge_{j=1}^T \overline{B_i}$, then the two graph sequences
$(\Gapp^1,\ldots,\Gapp^T)$ and $(\Gorg^1,\ldots,\Gorg^T)$ have the
same joint distribution. Since there is a fixed lower bound on the
probability for this event, any statement on
$(\Gapp^1,\ldots,\Gapp^T)$ which holds almost surely, also holds
almost surely for $(\Gorg^1,\ldots,\Gorg^T)$. The result follows.
\end{proof}
\begin{remark}
In the above claim we used the low probability for a step to get
omitted in order to show that statements that hold almost surely for
$\Gorg$ can be derived from such results on $\Gapp$. However, it is
worth noting that the two processes are much closer than that; if we
condition that $\Gapp^T$ has $M$ edges, then clearly it is
distributed as $\Gorg^M$. It is not difficult to show that, by the
low probability for omitting a step, $\Gapp(t)$ has the same
distribution as $\Gorg\left((1+o(1))t\right)$.
\end{remark}

\section{The behavior of the isolation ratio}\label{sec::isolation}
\subsection{Proof of Theorem \ref{I-thm}}
{We begin with the special case $K=0$, which we prove directly on
$\Gorg$. In this case, the differential equation \eqref{y-eq} takes
the simple form:
$$y'=-1~,~y(0)=1~,$$
and hence its unique solution is $y(t)=1-t$. Notice that, as $K$
equals $0$, the biased process connects two isolated vertices at
each step with probability $1$, as long as two such vertices exist.
Hence, at time $0 \leq t \leq 1$, there are $\lfloor t n/2 \rfloor$
edges which are vertex disjoint in pairs. Thus the ratio of isolated
vertices, $I(\Gorg(t)) =|\mathcal{C}_0(\Gorg(t))|/n$, satisfies:
$$|I(\Gorg(t))-\left(1 - t\right)| \leq \frac{2}{n}~,$$ where the $\frac{2}{n}$-term
is the rounding error.

We are left with the case $K > 0$. Fix $C > 0$, and let $G \sim
\Gapp(t)$ denote a graph at some point $t < C$ along the approximate
biased process $\Gapp$. We examine the effect that a single step of
$\Gapp$ has on the ratio isolated vertices,
$I(G)=|\mathcal{C}_0(G)|/n$.

Set $I=I(G)$; the total of the weights assigned to all ordered pairs
is $(I^2 + K(1-I^2))n^2$. Hence, with probability
$\frac{I(I-\frac{1}{n})}{I^2 + K(1-I^2)}$ the chosen edge is between
two formerly isolated vertices, and with probability
$\frac{2I(1-I)K}{I^2 + K(1-I^2)}$ precisely one end point of the
chosen edge was formerly isolated. Let $G'$ denote the graph after
performing the above step. The expected change in the isolation
ratio between $G$ and $G'$ satisfies:
\begin{eqnarray}\label{delta-I-eq}
 \mathbb{E}\left(I(G')-I(G)\right) &=& \left(-\frac{2}{n}\right) \frac{I(I-\frac{1}{n})}{I^2 + K(1-I^2)}
+\left(-\frac{1}{n}\right) \frac{2I(1-I)K}{I^2 + K(1-I^2)} =
\nonumber \\
&=& -\frac{2}{n} \left(1-\frac{(1-I)K}{I^2 + K(1-I^2)} -
\frac{I/n}{I^2 + K(1-I^2)}\right) ~.
\end{eqnarray}
Define: \begin{equation}\label{err-I-eq}\mathrm{err_y}(I,K) =
\frac{I/n}{I^2 + K(1-I^2)}~,
\end{equation}
and notice that the denominator of $\mathrm{err_y}(I,K)$ lies
between $1$ and $K$ for every value of $I$. Hence, as $K > 0$:
\begin{equation}\label{err-I-bound-eq}
  \mathrm{err_y}(I,K)\leq \frac{I/n}{\min\{1,K\}}=O(1/n)~.
\end{equation}
As $G'$ corresponds to $\Gapp(t+\frac{1}{n/2})$, we choose $\Delta
t = 2/n$ and rewrite \eqref{delta-I-eq} in the following form:
\begin{equation}\label{delta-exp-I-eq}
\mathbb{E}\frac{I(\Gapp(t+\Delta t))-I(\Gapp(t))}{\Delta t} =
\frac{(1-I)K}{I^2 + K(1-I^2)} - 1 + \mathrm{err_y}(I,K)~.
\end{equation} Notice that the left hand side of \eqref{delta-exp-I-eq}
resembles $\frac{d I}{d t}$, suggesting that the expected change in
the isolation ratio $I(\Gapp(t))$ is linked with the solution to the
differential equation \eqref{y-eq}. Notice that $y'$ is $C^\infty$,
and hence there is a unique solution to \eqref{y-eq}: we analyze
this solution in the next subsection, and now turn to show that
indeed it approximates $I(\Gapp(t))$. This will follow from a
general purpose theorem of \cite{WormaldDiffEq} (Theorem 5.1), which
we reformulate according to our needs, for the sake of simplicity (a
simpler version of the theorem appears in
\cite{WormaldDiffEqSimple}). If $\vec{y}=(y_1,\ldots,y_k)$, we use
the notation $(x,c\vec{y})$ to describe the tuple $(x,c y_1,\ldots,c
y_k)$.
\begin{theorem}[\cite{WormaldDiffEq}]\label{diff-method-thm}
Let $\vec{Y}=\vec{Y}(n)=(Y_1,\ldots,Y_l)$ denote $l\geq 1$ functions
from graphs on $n$ vertices to the real interval $[-C_0 n, C_0 n]$,
where $C_0 > 0$ is some constant, and let $\vec{f}=(f_1,\ldots,f_l)$
denote $l$ functions in $\mathbb{R}^{l+1}\rightarrow \mathbb{R}$.
Let $\{H^T\}$ denote a graph process on $n$ vertices beginning with
the edgeless graph. Let $\mathcal{D} \subset \mathbb{R}^{l+1}$ be a
bounded connected open set such that $(0,\frac{\vec{Y}(H^0)}{n})\in
\mathcal{D}$, and let $T_\mathcal{D}(H)$ denote the minimal time $T$
such that $(\frac{T}{n},\frac{\vec{Y}(H^T)}{n})$ no longer belongs to
$\mathcal{D}$. Assume the following:
\begin{enumerate}
\item (Boundedness Hypothesis) For some function $\beta=\beta(n)\geq
1$, $|\vec{Y}(H^{T+1})-\vec{Y}(H^T)|_\infty\leq \beta$ for every $T
\leq T_D(H)$.
\item (Trend Hypothesis) For some function
$\lambda_1=\lambda_1(n)=o(1)$, the following holds for every $T \leq
T_D(H)$: $\left|\mathbb{E}\left(\vec{Y}(H^{T+1})-\vec{Y}(H^T)\right)
- \vec{f}(\frac{T}{n},\frac{\vec{Y}(H^T)}{n})\right|_\infty \leq
\lambda_1$ .
\item (Lipschitz Hypothesis) There exists a
constant $L>0$ such that the following holds for every
$\vec{x},\vec{y}\in \mathcal{D}$:
$|\vec{f}(\vec{x})-\vec{f}(\vec{y})|_\infty \leq L
|\vec{x}-\vec{y}|_\infty$.
\end{enumerate}
Then the following holds: \begin{enumerate}\item
\label{diff-method-result1}There exists a unique solution to the
following system of $l$ differential equations:
$$\left\{\begin{array}{l}
\frac{d u_1}{d x}=f_1(x,\vec{u}),~\ldots,~\frac{d
u_l}{d x}=f_l(x,\vec{u})\\
\vec{u}(0)=\frac{\vec{Y}(H^0)}{n}\\
\end{array} \right. ~,$$
where $\vec{u}$ denotes $(u_1,\ldots,u_l)$. Let
$\tilde{u}_1(x),\ldots,\tilde{u}_l(x)$ denote this solution.
\item \label{diff-method-result2} Let $\lambda > \lambda_1$,
$\lambda=o(1)$, and let $\sigma > 0$ be such that
$(x,\vec{\tilde{u}}(x))$ is at least $C_1 \lambda$ away from the
boundary of $\mathcal{D}$ for every $x \in [0,\sigma]$, where
$C_1>0$ is some constant. Then there exists a constant $C_2>0$, such
that with probability
$1-O(\frac{\beta}{\lambda}\exp(-\frac{n\lambda^3}{\beta^3}))$, the
following holds for every $0\leq T \leq \sigma n$ and every $1 \leq
i \leq l$: $|\frac{Y_i(H^T)}{n} - \tilde{u}_i(\frac{T}{n})| \leq C_2
\lambda$.
\end{enumerate}
\end{theorem}
\begin{remark}
Note that for the purpose of proving Theorem \ref{I-thm} we will only need the special case
$l=1$ of Theorem \ref{diff-method-thm} above. However, in Section
\ref{sec::susceptibility}, when analyzing the
critical point for the emerging time of the giant component,
we apply the above theorem for $l=2$.
\end{remark}
We claim that the following substitution completes the proof of
Theorem \ref{I-thm}:
\begin{eqnarray} Y(G)&=&|\mathcal{C}_0(G)|~,
\nonumber \\
f(x,y)&=&2(\frac{(1-y)K}{y^2+(1-y^2)K}-1) ~.\nonumber
\end{eqnarray} For the set $\mathcal{D}$ we choose a bounded connected open set
containing the rectangle $[0,C]\times[0,1]$. Notice that $0 \leq
\frac{Y(G)}{n}\leq 1$ for every graph $G$ on $n$ vertices, hence
$(\frac{T}{n},\frac{Y(\Gapp^T)}{n})$ belongs to $\mathcal{D}$ for
every $1 \leq T \leq C n/2$ and $T_D(\Gapp)$ is at least $C$.

Indeed, $|Y(\Gapp^{T+1})-Y(\Gapp^T)| \leq 2$, thus a choice of
$\beta=2$ confirms the Boundedness Hypothesis. Next,
\eqref{delta-I-eq}, \eqref{err-I-eq} and \eqref{err-I-bound-eq}
imply that
$$|\mathbb{E}\left(Y(\Gapp^{T+1})-Y(\Gapp^T)\right)-f(\frac{T}{n},\frac{Y(\Gapp^T)}{n})|
= 2 \mathrm{err_y}(\frac{Y(\Gapp^T)}{n},K) = O(1/n)~.$$ Thus,
setting $\lambda_1=\frac{1}{n}$ verifies the Trend Hypothesis.
Finally, $f(x,y)$ is clearly $C^\infty$ (recall that $K>0$) and
hence satisfies the Lipschitz condition.

Proposition \ref{prop-y-sol-exp-bound}, proved in the next
subsection by analyzing the differential equation \eqref{y-eq},
states that the solution $\tilde{u}(x)$ is bounded between
$\exp(-x/Q)$ and $\exp(-Q x)$. Hence, for every $0 \leq x \leq C $
the solution $\tilde{u}(x)$ is at least a constant away from the
boundary of $\mathcal{D}$, and we can easily choose $\sigma=C$
(regardless of our choice of $\lambda$).

Altogether, a choice of $\lambda = n^{-1/4}$ implies that with
probability $1-O(n^{1/4} \exp(-n^{1/4}))=1-o(\exp(-\frac12 n^{1/4})$,
the quantity $Y(G)$ satisfies:
$$ \frac{Y(\Gapp^T)}{n} = \tilde{u}\big(\frac{T}{n}\big)+O(n^{-1/4})~,$$
for every $0 \leq T \leq C n / 2$, where $\tilde{u}$ is the (unique)
solution to the equation $ \frac{d u}{d x} = f(x, u)$.

Scaling the time to units of $n/2$ edges, we obtain that, with
probability $1-o(\exp(-\frac12 n^{1/4})$, the following holds for every
$0 \leq t \leq C$:
$$ |I(\Gapp(t)) - y(t)|= O(n^{-1/4})~,$$
where $y(x)=\tilde{u}(x/2)$. Hence, $y$ is the unique solution to
the equation: $$\frac{d y}{d x} = \frac{1}{2}\frac{d u}{d x} =
\frac{(1-y)K}{y^2+(1-y^2)K}-1~,$$ with the starting condition
$y(0)=\tilde{u}(0)=1$, completing the proof.
 \qed}

\subsection{The behavior of $y(t)$ when $K>0$}
The behavior of $y(t)$ along the biased process is crucial to the
understanding of how the susceptibility grows, as we show in the
next section. The behavior of the ratio of isolated vertices (and
the corresponding function $y(t)$) was already stated for the case
$K=0$. We thus assume $K>0$, and set $Q = \lceil
\max\{K,\frac{1}{K}\}\rceil$. By Corollary
\ref{cor-stochastic-domination}, we obtain that for every $\epsilon
> 0$ and every $C > 0$, the following holds with high probability for every $t \in [0,C]$:
\begin{equation}
\label{exponential-I-bounds}  \mathrm{e}^{-Q t} - \epsilon \leq
I(\Gorg(t)) \leq \mathrm{e}^{-t / Q} + \epsilon~.
\end{equation} This follows from the uniform continuity of the
functions $\exp(-t/Q)$ and $\exp(-Q t)$, which describe
$I(\Gorg[1](t/Q))$ and $I(\Gorg[1](Q t))$ respectively.

We next show that a stronger result than
\eqref{exponential-I-bounds} can be easily derived directly from the
differential analysis of $y(t)$, as stated by the next proposition:
\begin{proposition}\label{prop-y-sol-exp-bound} The ODE
\eqref{y-eq} for $K>0$ has the following properties:
\begin{enumerate}
\item Its solutions are strictly monotone increasing in $K$. That
is, if $y_1(t), y_2(t)$ are the solutions that correspond to $K_1 <
K_2$, then $y_1(t) < y_2(t)$ for all $t > 0$.
\item The solution $y(t)$ is strictly monotone decreasing in $t$,
and satisfies:
\begin{equation}\label{y-sol-exp-bound}
\exp\left(-\frac{t}{\min\{1,K\}}\right) \leq y(t) \leq
\exp\left(-\frac{t}{\max\{1,K\}}\right)~,
\end{equation}
for all $t \geq 0$. Furthermore, the inequalities are strict
whenever $K \neq 1$.
\end{enumerate}\end{proposition}
\begin{proof}
We begin by proving an inequality analogous to
\eqref{exponential-I-bounds}, stating that whenever $K \neq 1$,
$y(t)$ is strictly monotone decreasing in $t$ and satisfies:
\begin{equation}\label{eq-initial-y-sol-exp-bound}
\exp(-Q t) < y(t) < \exp(-t/Q)
\end{equation}
for any $t > 0$, where $Q = \max\{\frac{1}{K},K\}$. Inequality
\eqref{y-sol-exp-bound} will follow directly from
\eqref{eq-initial-y-sol-exp-bound} once we prove that $y$ is
strictly monotone increasing in $K$, since $\exp(-t)$ is the
solution to the ODE \eqref{y-eq} for $K=1$.

Suppose $K>1$ and let $u(t) = \exp\{-t/Q \}$. First, we prove that
$u$ satisfies the inequality
\begin{equation}\label{ineq}
u' > \frac{(1-u)K}{u^2 + (1-u^2)K}-1,\quad \mbox{for } t>0.
\end{equation}
Indeed, $u' = -\ds\frac{1}{Q}\exp\{-t/Q \} = -\ds\frac{u}{Q}$.
Therefore, (\ref{ineq}) holds iff
\begin{equation}\label{equiv-ineq}
-\frac{u}{Q} > \frac{(1-u)K}{u^2 + (1-u^2)K}-1, \quad \mbox{for }
0<u<1.
\end{equation}
which, after some manipulations, is equivalent to the inequality
\[
p_Q(u)\equiv (K-1)u^2 - Q(K-1)u + K(Q-1) > 0, \quad \mbox{for }
0<u<1.
\]
The quadratic polynomial $p_Q(u)$ satisfies $p_Q(u)\to+\infty$ as
$u\to \pm\infty$, because $K>1$. Furthermore, if $Q>1$ then $p_Q(0)
= K(Q-1)
>0$ and $p_Q(1)=Q-1>0$. Therefore, if $p_Q(u)$ does not have a root in
the interval $(0,1)$ then the inequality holds. In the special case
of $Q=K$ we have $p_K(u) = (K-1)\left(u^2 - Ku + K \right)$. The
roots of $p_K(u)$ are $u_\pm  = \ds\frac{K}{2}\pm
\ds\frac{1}{2}\sqrt{K^2-4K}$. For $1<K<4$ there are no real roots.
For $K=4$ there is a double root at $u=2$ and for $K>4$ there are
two distinct roots $u_\pm > 1$. Therefore, the inequality
(\ref{ineq}) holds for all values of $Q=K>1$ as asserted. By the
equivalent inequality \eqref{equiv-ineq} (recall that $0<y<1$ for
$t>0$), we obtain that $y'(t)<0$ for every $t\geq 0$, and hence
$y(t)$ is strictly monotone decreasing.

To complete the proof of the upper bound of
\eqref{eq-initial-y-sol-exp-bound}, let $f(w) = \ds\frac{(1-w)K}{w^2
+ (1-w^2)K}-1$. The differential equation (\ref{y-eq}) and
inequality (\ref{ineq}) imply that the functions $y(t)$ and $u(t)$
satisfy for $t>0$
\begin{eqnarray}
y' &=& f(y), \nonumber \\
u' &>& f(u), \nonumber
\end{eqnarray}
together with the mutual initial condition $y(0)=u(0)=1$. By
standard analytical considerations, this implies that $u(t)>y(t)$
for all $t>0$ (to see this, note that $u(t)>y(t)$ for $t\in
[0,\epsilon]$ and some small $\epsilon>0$, and that the smallest
point $t^*>\epsilon$ satisfying $u(t^*)=y(t^*)$ cannot satisfy
$u'(t^*)-y'(t^*) > f(u(t^*))-f(y(t^*))=0$, yielding a
contradiction).

The lower bound in \eqref{eq-initial-y-sol-exp-bound} and the case
$0<K<1$ are similar.

To prove that $y$ is strictly monotone increasing in $K$, take
$0<K_1<K_2$, and let $y_1$ and $y_2$ the solutions of equation
\eqref{y-eq} for $K_1$ and $K_2$ respectively. Considering $f$,
defined as above, as a function of both $w$ and $K$, $f$ satisfies:
$${\frac{\partial f}{\partial
K}(w,K)=\frac{(1-w)w^2}{(w^2+(1-w^2)K)^2}}~,$$ hence $\frac{\partial
f}{\partial K}(w,K)>0$ for every $0<w<1$. Recall that inequality
\eqref{eq-initial-y-sol-exp-bound} guarantees that $y(t)>0$ for
every $t$. Therefore, for any $t>0$, $f(y(t),K)$ is monotone
increasing in $K$, and in particular, for any $t>0$ we have
$y_1'=g(y_1)$ and $y_2'> g(y_2)$, where $g(y)=f(y,K)$. Thus, the
above argument for proving inequality
\eqref{eq-initial-y-sol-exp-bound} completes the proof.
\end{proof}

We note that the bounds of \eqref{y-sol-exp-bound} are much weaker
than the estimations which can be obtained by examining the
asymptotic behavior of the ODE \eqref{y-eq}, as we proceed to do in
Section \ref{sec::critical}.

\section{The susceptibility of the biased
process}\label{sec::susceptibility}
\subsection{Proof of Proposition \ref{z-thm}}
We prove Proposition \ref{z-thm} by showing that there exists an
$M>0$ such that $0<t_c(K)<M$ for every $K>0$. The fact that $z(t)$
is continuous on $[0,t_c)$ follows from standard considerations in
differential analysis. The case $K=1$ is trivial, as in this case
equation $\eqref{z-eq}$ takes the simple form
$$z'=z^2~,~z(0)=1~,$$
hence its solution is $1/(1-t)$ and $t_c(1)=1$.

Take $K>1$, and let $z(t)$ denote the (unique) solution to equation
\eqref{z-eq}. Recall that by Proposition \ref{prop-y-sol-exp-bound},
$y(t)$, the solution to equation \eqref{y-eq}, satisfies:
$$
\mathrm{e}^{-t} < y(t) < \mathrm{e}^{-t / K}
$$
for every $t > 0$. Thus, $y(t)>0$ for every $t\geq 0$, and we obtain
that $\ds{\frac{K}{y^2+(1-y^2)K}>1}$ for every $t\geq 0$. Therefore:
$$z' > z^2 ~,~ z(0)=1~,$$
and by the method used in the proof of Proposition
\ref{prop-y-sol-exp-bound}, we obtain that $z(t) \geq 1/(1-t)$ for
every $t\geq 0$. We deduce that $0<t_c(K) \leq 1$ for every $K > 1$.

Let $0 < K < 1$, and again let $y(t)$ and $z(t)$ denote the
solutions to equations \eqref{y-eq} and \eqref{z-eq} respectively.
Recalling that $y(0)=1$ and $y(t)$ decreases to $0$ as $t\to\infty$,
let $t^*>0$ be such that $y(t^*)=\sqrt{K}$. It is easy to verify
that the derivative of $f(y)=\frac{(1-y)K}{y^2+(1-y^2)K}-1$ is
strictly negative for every $0\leq y \leq 1$. Thus, \eqref{y-eq}
implies that for every $0 \leq t \leq t^*$,
\begin{equation}
  y'(t) = f(y) \leq f(\sqrt{K}) = \frac{1-\sqrt{K}}{2-K}~.
\end{equation}
By the Mean Value Theorem, we obtain that:
\begin{equation}\label{t^*-eq}
  t^* \leq
  \frac{1-\sqrt{K}}{1-\frac{1-\sqrt{K}}{2-K}}=\frac{1}{\frac{1}{1-\sqrt{K}}-\frac{1}{2-K}} \leq
  2~,
\end{equation}
where the last inequality is by the fact that the function
$\frac{1}{1-\sqrt{K}}-\frac{1}{2-K}$ is monotone increasing from
$\frac{1}{2}$ to $\infty$ as $K$ goes from $0$ to $1$. In addition,
by \eqref{z-eq}, $z'(t)>0$ for every $t\geq 0$
(this applies to every $K\geq 0$). Therefore:
\begin{equation}\label{hat-bound-1}
  z(t^*) > z(0) = 1~.
\end{equation}
Hence: \begin{equation}\label{hat-bound-2} c(t) =
\frac{K}{y^2+(1-y^2)K} > \frac{1}{2}
\end{equation}
for every $t \geq t^*$. By \eqref{hat-bound-1} and
\eqref{hat-bound-2}, defining $\hat{z}(t) = z(t+t^*)$ gives:
\begin{equation}
\hat{z}' > \frac{1}{2} \hat{z}^2 ~,~ \hat{z}(0)>1~,
\end{equation}
and by comparing $\hat{z}$ with the function $\frac{2}{2-t}$ we
deduce that $\hat{z}$ has a singularity point $\hat{t}_c \leq 2$.
Altogether, $z(t)$ satisfies $t_c \leq \hat{t}_c+t^* \leq 4$,
completing the proof. \qed

\begin{remark} Theorem \ref{critical-thm} states that
$t_c(K)$ is in fact monotone decreasing as a function of $K$ and
thus bounded by $t_c(0)=\frac{3}{2}$. We note that Theorem
\ref{S-thm}, used in the proof of Theorem \ref{critical-thm},
requires a bound on $t_c(K)$.\end{remark}

\subsection{The susceptibility at the subcritical phase}
In this subsection we combine Theorem \ref{thm-M-bounded} with
methods from \cite{SpencerWormald} to prove part \ref{S-thm-sub} of
Theorem \ref{S-thm}. The proof relies on the fact that a bounded
susceptibility ensures a logarithmic upper bound for the components,
and this in turn ensures that the solution to the differential
equation \eqref{z-eq} stays a good approximation for the
susceptibility. Hence, the fact that $S(\Gorg(t))$ is approximated
by $z(t)$ and the fact that every component is of size $O(\ln n)$
ensure each other along the subcritical phase.

As we next show, the special case $K=0$ again follows from the behavior of $S(\Gorg[1])$.
Indeed, \eqref{z-eq} takes the following simple form when $K=0$:
$$z'=1~,~z(0)=1~,$$
and hence its unique solution is $z(t)=1+z$. By definition, as
$K=0$, as long as there exists a pair of isolated vertices in the
graph, such a pair forms the next edge. This amounts to trading two
components of size $1$ with one component of size $2$, and hence at
each such step the susceptibility, being the sum of squares of the
component sizes, increases by $\frac{2}{n}$. Thus, at time $0 \leq t
\leq 1$, the susceptibility equals $1+\frac{2}{n} \lfloor t n/2
\rfloor$, hence:
$$|S(\Gorg(t))-\left(1 + t\right)| \leq \frac{2}{n}~,$$ where the $\frac{2}{n}$-term
is the rounding error. As of this point, the process is equivalent
to $\Gorg[1]$ on $\lfloor \frac{n+1}{2} \rfloor$ vertices, $\lfloor
n/2 \rfloor$ of which represent components of size $2$. Hence, the
result on $S(G)$ from this point on is derived from the case $K=1$
of Theorem \ref{S-thm}. The solution for the differential equation
\eqref{z-eq} when $K$ equals $1$ is $z=\frac{1}{1-t}$, hence, we
obtain that the susceptibility is within $o(1)$ of the function
$\hat{z}=\frac{2}{1-2(t-1)}=\frac{1}{3/2-t}$ for $t \geq 1$.

Throughout the remainder of the proof, assume therefore that $K>0$.
We need the following definition: a graph $G$ is said to have a {\em
$k,c$ component tail} if, for every $s$, the probability that a
uniformly chosen vertex has a component of size at least $s$ is at
most $k\exp(-c s)$:
$$\Pr_{v \in V(G)}[|\mathcal{C}_v|\geq s] \leq k \mathrm{e}^{-c
s}~.$$ Notice that if a graph $G$ on $n$ vertices has this property
then its largest component is of size at most $\max\{1,1/c\} \ln n$,
provided that $n$ is sufficiently large.

The following theorem (\cite{SpencerWormald}, Theorem 3.1) is
crucial to the proof:
\begin{theorem}[\cite{SpencerWormald}]\label{thm-k-c}
Let $\ell,k,c$ be positive real numbers. Let $G$ be a graph on $n$
vertices with a $k,c$ component tail. Let $H$ be a random graph with
edge probability $p = \frac{\hat{t}}{n}$ on the same vertex set,
where $\hat{t}$ is fixed. Set $G^+ = G \cup H$. \begin{enumerate}
\item (Subcritical) Assume $S(G) \leq \ell$ for all $n$. Let $\hat{t} \ell
< 1$. Then there exist $k^+$, $c^+$ (dependent on $k,c,\ell,\hat{t}$
but not on $n$ nor $G$) such that $G^+$ almost surely has a
$k^+$,$c^+$ component tail. In particular, all components have size
$O(\ln n)$.
\item (Supercritical) Assume $S(G) > \ell$. Let $\hat{t} \ell > 1$. Then
$G^+$ almost surely has a giant component. More precisely, there
exists $\gamma
> 0$ (dependent on $k,c,\ell,\hat{t}$ but not on $n$ nor $G$) such that $G^+$ has a
component of size at least $\gamma n$.\end{enumerate}\end{theorem}

The following lemma proves part \ref{S-thm-sub} of Theorem
\ref{S-thm} (analogous to Theorem 4.2 of \cite{SpencerWormald}):
\begin{lemma}\label{lemma-subcritical}
Let $\epsilon > 0$. For every $t \in [0,t_c-\epsilon]$, where $t_c$
is as defined in Proposition \ref{z-thm}, the following two
conditions hold almost surely:
\begin{enumerate}
  \item \label{sub-lemma-susc-app} (Susceptibility Approximation)
  $|S(\Gorg(t)) - z(t)| = o(1)$, where $z(t)$ is the solution to the
  differential equation \eqref{z-eq}.
  \item \label{sub-lemma-small-comp} (Small Components) $\Gorg(t)$
  has a $k,c$ component tail, for some $k=k(t),c=c(t)$.
\end{enumerate}
\end{lemma}
\begin{proof}
 Set $Q = \lceil \max\{\frac{1}{K},K\}\rceil$, and let:
 \begin{equation}
   L=\lceil (z(t_c-\epsilon)+\epsilon)Q t_c\rceil~,~t_j =
   \frac{t_c-\epsilon}{L}j~(j=0,\ldots,L)~.
 \end{equation}

Assume by induction that properties \ref{sub-lemma-susc-app} and
\ref{sub-lemma-small-comp} of the lemma almost surely hold for $t
\in [0,t_j]$ for some $0 \leq j < L$ (the case $j=0$ is trivial). We
show that with probability $1-o(1)$ properties
\ref{sub-lemma-susc-app} and \ref{sub-lemma-small-comp} hold for the
interval $[t_j,t_{j+1}]$ as well, and the lemma follows from a union
bound on the complement events.

First, consider the Small Components property. By the induction
hypothesis and the monotonicity of $z(t)$ (recall that
$z'(t)>0$ for all $t\geq 0$,
as mentioned in the proof of Proposition \ref{z-thm}),
 $$S(\Gorg(t_j)) \leq z(t_j) + o(1) \leq z(t_L) +
o(1)~.$$ Let: \begin{equation}\label{S-bound-of-tQ} \ell =
z(t_L)+\epsilon~,~\hat{t} = Q \frac{t_L}{L}~. \end{equation} With
high probability, $\Gorg(t_j)$ has a $k,c$ component tail for some
constants $k,c$ (dependent on $t_j$), $S(\Gorg(t_j))\leq \ell$ and
$\hat{t}\ell \leq \frac{t_c-\epsilon}{t_c} < 1$. Thus, Theorem
\ref{thm-k-c} implies that $\Gorg[1](\hat{t})|_{\Gorg(t_j)}$ almost
surely has a $k^+,c^+$ component tail for some constants $k^+,c^+$.
By Theorem \ref{thm-M-bounded} (inequality \eqref{dom-M-by-1}) we
obtain that $\Gorg(\hat{t}/Q)|_{\Gorg(t_j)}$ almost surely has a
$k^+,c^+$ component tail, hence $\Gorg(t_{j+1})$ satisfies the Small
Components property. Thus, in particular, with high probability all
the components of $\Gorg(t_{j+1})$ are of size $O(\ln n)$.

We next analyze the change in the susceptibility along the
approximate biased process. Since $\Gorg^T$ is equivalent to
$\Gapp^{T'}$ for some $T'>T$ (where $T'-T$ is the number of
redundant steps in $\Gapp$), the Small Components property which we
proved above implies that, with high probability, the largest
component of $\Gapp(t_{j+1})$ is of size at most $\alpha \ln n$ for
some $\alpha > 0$. Assume therefore that this indeed is the case.

Let $G$ denote an instance of $\Gapp(t)$ on $n$ vertices for some
$t<t_{j+1}$, and set $I=I(G)$ and $S=S(G)$. Recall that the total of
the weights assigned to all ordered pairs is $(I^2 + K(1-I^2))n^2$,
and that in case that the chosen pair in the next step is a loop or
an edge which already exists in $G$, the step is omitted. Let $G'$
denote the graph after the next step, and let $A_{C_1,C_2}$ denote
the event in which the newly chosen pair joins the two components
$C_1,C_2 \in \mathcal{C}(G)$. If $G'$ equals $G$ or if the new edge
is an internal edge of some component $C \in \mathcal{C}(G)$ (the
event $A_{C,C}$ occurred), then $S(G')=S(G)$. Otherwise, a new edge
between $C_1,C_2\in\mathcal{C}(G)$ implies:
\begin{equation}\label{susceptibility-delta-eq}S(G')
= \frac{1}{n}\sum_{C\in \mathcal{C}(G')}|C|^2 = S + \frac{2}{n}|C_1|
|C_2|~.
\end{equation} Hence: \begin{equation}\label{delta-exp-S-initial-eq}
\frac{\mathbb{E}(S(G')-S(G))}{2/n}=\sum_{C_1\in\mathcal{C}}
\mathop{\sum_{C_2\in\mathcal{C}}}_{C_1\neq C_2}
|C_1||C_2|\Pr[A_{C_1,C_2}]~.\end{equation} There are two cases to
consider, according to which we divide the summands of the right
hand side of \eqref{delta-exp-S-initial-eq} to two terms, $\Delta_1$
and $\Delta_2$. In the first case, a new edge $(u,v)$ is added
between two distinct and formerly isolated vertices $u,v$. Since
$|\mathcal{C}_u|=|\mathcal{C}_v|=1$, the contribution of this case
to the sum in \eqref{delta-exp-S-initial-eq} is its probability:
\begin{equation}\label{exp-S-delta1-eq}\Delta_1=\sum_{u \neq v \in
\mathcal{C}_0}\Pr[A_{\mathcal{C}_u,\mathcal{C}_v}]=\frac{I^2}{I^2+(1-I^2)K}
- \Pr[\cup_{u \in \mathcal{C}_0}
A_{\mathcal{C}_u,\mathcal{C}_u}]~.\end{equation} In remaining case,
$A_{C_1,C_2}$ occurs where $C_1\neq C_2$ and at least one of the
components $C_1,C_2$ is not an isolated vertex. The contribution to
the right hand side of \eqref{delta-exp-S-initial-eq} in this case
is:
\begin{equation}\label{exp-S-delta2-eq}
\Delta_2= \mathop{\sum_{C_1 \in \mathcal{C}}}_{C_2 \in
\mathcal{C}}\frac{K|C_1|^2|C_2|^2}{(I^2+(1-I^2)K)n^2}-\frac{I^2
K}{I^2+(1-I^2)K} - \mathop{\sum_{C \in
\mathcal{C}}}_{|C|>1}\frac{K|C|^4}{(I^2+(1-I^2)K)n^2}.
\end{equation}
Define:
\begin{eqnarray}
  \mathrm{err_z}(I,K)&=&\Pr[\cup_{u \in \mathcal{C}_0}
A_{\mathcal{C}_u,\mathcal{C}_u}] + \mathop{\sum_{C \in
\mathcal{C}}}_{|C|>1}\frac{K|C|^4}{(I^2+(1-I^2)K)n^2}= \nonumber \\
&=& \frac{K}{I^2+(1-I^2)K} \sum_{C \in \mathcal{C}}
\frac{|C|^4}{n^2} +
\frac{1-K}{I^2+(1-I^2)K}\frac{I}{n}~.\label{err-S-eq}
\end{eqnarray}
Combining \eqref{exp-S-delta1-eq},\eqref{exp-S-delta2-eq} and
\eqref{err-S-eq}, the following holds:
$$ \Delta_1 + \Delta_2 =
\frac{(1-K)I^2}{I^2+(1-I^2)K} + \frac{K}{I^2+(1-I^2)K}
\mathop{\sum_{C_1 \in \mathcal{C}}}_{C_2 \in \mathcal{C}}
\frac{|C_1|^2|C_2|^2}{n^2} - \mathrm{err_z}(I,K)~,$$ and therefore:
\begin{equation}\label{delta-exp-S-eq}
\frac{\mathbb{E}(S(G')-S(G))}{2/n}= \frac{(1-K)I^2}{I^2+(1-I^2)K} +
\frac{K}{I^2+(1-I^2)K} S^2(G) - \mathrm{err_z}(I,K)~.
\end{equation}
This explains the choice of the differential equation \eqref{z-eq},
provided that $\mathrm{err_z}$ is redundant; this is ensured by the
fact that the largest component of $G$ is of size at most $\alpha \ln
n$. To see this, notice that $G$ satisfies $\frac{1}{n}\sum_{C \in
\mathcal{C}(G)}|C|^4 \leq \sum_{C \in
\mathcal{C}(G)}\frac{|C|}{n}(\alpha \ln n)^3 = (\alpha \ln n)^3$. We
obtain the following bound on $\mathrm{err_z}$:
\begin{equation}\label{err-S-bound-eq}
\mathrm{err_z}(I,K) \leq \frac{K}{(I^2+(1-I^2)K} \frac{\alpha^3
\ln^3 n}{n} + \frac{1-K}{I^2+(1-I^2)K}\frac{I}{n} = O(\frac{\ln^3
n}{n})~.
\end{equation}

It is left to verify the conditions of Theorem
\ref{diff-method-thm}. Although we may assume that the statements of
Theorem \ref{I-thm} hold with the respect to the approximate biased
process, and $|I(\Gapp(t))-y(t)|=O(n^{-1/4})$ for every $0 \leq t
\leq t_c$, where $y$ is the unique solution to \eqref{y-eq}, the
$O(n^{-1/4})$ approximation error will not be sufficient for proving
that $S(G)$ is approximated by $z(G)$. The difficulty is in proving
the Trend Hypothesis; $d S(G)/d t$ is approximated by
$g_1(I)+g_2(I)S^2(G)$ for some continuous functions $g_1,g_2$, and
using the $y$-approximation for $I(G)$ implies an $O(n^{-1/4})
S^2(G)$ error. Therefore, we apply Theorem \ref{diff-method-thm} on
both $I(G)$ and $S(G)$ (re-proving the result on $I(G)$ for the
interval $[0,t_{j+1}]$).

Keeping the notations of the theorem, define:
\begin{eqnarray} Y_1(G)&=& n I(G) = |\mathcal{C}_0(G)|~,
\nonumber \\
Y_2(G)&=&n S(G) = \sum_{C \in \mathcal{C}(G)}|C|^2~,
\nonumber \\
f_1(x,y,z)&=&2(\frac{(1-y)K}{y^2+(1-y^2)K}-1) ~,\nonumber \\
f_2(x,y,z)&=&2(\frac{K}{y^2+(1-y^2)K}\left(z^2-1\right)+1)
~.\nonumber
\end{eqnarray} The set $\mathcal{D}$ will be a bounded connected open
set containing the domain $[0,t_{j+1}]\times[0,1]\times[0,\ell]$,
where $\ell$ is the constant defined in \eqref{S-bound-of-tQ}.
Clearly, $\mathcal{D}$ contains
$(\frac{T}{n},\frac{\vec{Y}(\Gapp^T)}{n})$ as long as $T \leq
t_{j+1} n/2$ and $S(\Gapp^T) \leq \ell$, and in particular, for
$T=0$.

Recall that the first coordinate of $\vec{Y}$ satisfies the
Boundedness Hypothesis for a choice of $\beta=2$ and every $T \leq
t_{j+1} n/2$. Since every component of $G$ is of size at most
$\alpha \ln n$, \eqref{susceptibility-delta-eq} implies that
$|Y_2(\Gapp^{T+1})-Y_2(\Gapp^T)| \leq 2(\alpha \ln n)^2$, thus a
choice of $\beta=2\alpha^2 \ln n$ verifies the Boundedness
Hypothesis for every $T \leq t_{j+1} n/2$.

The Trend Hypothesis was satisfied by $Y_1,f_1$ in the proof of
Theorem \ref{I-thm} using a choice of $\lambda_1=1/n$ ($f_1$ is
independent of $z$). By \eqref{delta-exp-S-eq} and
\eqref{err-S-bound-eq}, we have:
$$ |\mathbb{E}\left(Y_2(\Gapp^{T+1})-Y_2(\Gapp^T)\right)-f_2(\frac{T}{n},\frac{\vec{Y}(\Gapp^T)}{n})|
= 2 \mathrm{err_z}(I,K) = O(\frac{\ln^3 n}{n})~.$$ Therefore, a
choice of $\lambda_1=O(\frac{\ln^3 n}{n})$ verifies the Trend
Hypothesis.

The Lipschitz condition is again satisfied by the fact that
$f_i(x,y,z)$, $i=1,2$, are clearly $C^\infty$ for $K>0$.

Notice that our definition of $\ell$, relying on Proposition
\ref{z-thm}, is such that the unique solution to the system of
differential equations remains at least $\epsilon$-far from the
boundary of the domain $\mathcal{D}$ for every $0 \leq x \leq
t_{j+1}$, hence a choice of $\sigma=t_{j+1}$ is valid.

Altogether, a choice of $\lambda = n^{-1/4}$ implies that with
probability $1-o(\exp(-\frac12 n^{1/4}))$ there exists some constant
$C_2>0$ such that the following holds for every $0 \leq T \leq
t_{j+1} n/2$:
$$ \left|\frac{\vec{Y}(\Gapp^T)}{n} -
(\tilde{y}(\frac{T}{n}),\tilde{z}(\frac{T}{n}))\right|_\infty \leq
C_2 n^{-1/4}~,$$ where $(\tilde{y},\tilde{z})$ is the (unique)
solution to the equation system:
$$\left\{\begin{array}{l}
\frac{d y}{d x}=f_1(x,y,z),~\frac{d
z}{d x}=f_2(x,y,z)\\
y(0)=z(0)=1\\
\end{array} \right. ~.$$
Scaling the time to units of $n/2$ edges, we obtain that, with
probability $1-o(\exp(-\frac12 n^{1/4}))$, the following holds for every
$0 \leq t \leq t_{j+1}$:
$$ |S(\Gapp(t)) - z^*(t)|= O(n^{-1/4})~,$$
where $z^*(x)=\tilde{z}(x/2)$. Hence, $z^*$ is the unique solution
to the equation: $$\frac{d z^*}{d x} = \frac{1}{2}\frac{d \hat{z}}{d
x} = \frac{K}{\tilde{y}^2+(1-\tilde{y}^2)K}z^2+
\frac{(1-K)\tilde{y}^2}{\tilde{y}^2+(1-\tilde{y}^2)K}~,$$ with the
initial condition $z^*(0)=\tilde{z}(0)=1$, completing the proof of
the lemma.
\end{proof}

\subsection{The susceptibility at the supercritical phase}
\begin{proof}[\textbf{\em Proof of Theorem \ref{S-thm} part \ref{S-thm-super}}]
Let $Q =\lceil \max\{K,\frac{1}{K}\}\rceil$, and assume without loss
of generality that $\epsilon < t_c$. By Proposition \ref{z-thm}
there exists some $0 < t^- < t_c$ satisfying $z(t^-) =
\frac{Q+2}{\epsilon}$. Lemma \ref{lemma-subcritical} implies that,
with high probability,
$$S(\Gorg(t^-)) \geq \frac{Q+2}{\epsilon}(1-o(1)) > \frac{Q+1}{\epsilon} ~,$$
and in addition, $\Gorg(t^-)$ has a $k,c$ component tail. Fix
$H=\Gorg(t^-)$. Setting $\ell = \frac{Q+1}{\epsilon}$ and
$\hat{t}=\frac{\epsilon}{Q}$, Theorem \ref{thm-k-c} implies that
$\Gorg[1](\frac{\epsilon}{Q})|_H$ almost surely has a giant
component. Hence, by Theorem \ref{thm-M-bounded} (inequality
\eqref{dom-1-by-M}), $\Gorg(\epsilon)|_H=\Gorg(t^-+\epsilon)$ almost
surely has a giant component, and in particular,
$\Gorg(t_c+\epsilon)$ has a giant component with high probability.

\end{proof}

\subsection{Computer simulations and numerical results}
We have implemented an efficient simulator of the model $\Gorg$,
which runs in a time complexity of $O(n \log n)$, and uses a data
structure of size $O(n)$. Maintaining the sets of isolated and
non-isolated vertices allows an $O(1)$ cost for randomizing the next
edge at each step whenever one of its endpoints is isolated; the
cost of selecting an edge between two non-isolated vertices is
$O(\log n)$, where the set of existing edges is kept in a balanced
tree, which is an efficient data-structure. Finally, maintaining the
connected components in linked lists according to the Weighted-Union
Heuristic (see, e.g., \cite{CLR} p. 445) provides an average cost of
$O(\log n)$ for uniting components.

We conducted a series of simulations with $n=10^6$, and values of
$K$ ranging from $0$ to $10$ in steps of $0.1$. In order to assess
the critical time $t_c(K)$, we considered the minimal time at which
the largest component of $\Gorg$ was of size $\alpha n$, where
$\alpha=0.01$. These estimations of $t_c(K)$ were supported by tests
of higher values of $n$, such as $2.5\cdot10^6$ and $5\cdot10^6$
vertices. We note that the constant $\alpha$ decreases with $K$
(independently of $n$). Therefore, automatically assessing $t_c(K)$
for large values of $K$ requires a different method, such as a
threshold for the derivative of the size of the largest component,
or for the derivative of the susceptibility.

\begin{figure}
\centering
\includegraphics{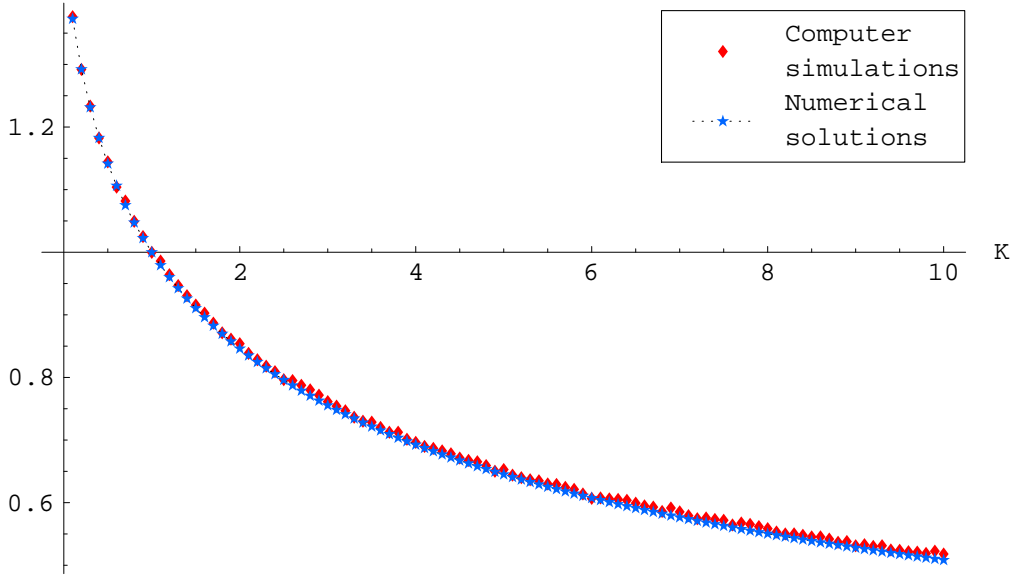}
\caption{Comparison of the numerical results for $t_c(K)$, and
estimations of $t_c(K)$ according to computer simulations of the
model.} \label{fig::sim}
\end{figure}

The simulated value of $t_c(K)$ in each setting was averaged over
$10$ tests conducted as mentioned above. Figure \ref{fig::sim} shows
an excellent agreement between these values and the numerical
approximations of $t_c(K)$, obtained by numerical solutions of the
ODEs \eqref{y-eq} and \eqref{z-eq} by Mathematica.

\section{Analytical properties of the critical time}
\label{sec::critical} In this section we prove Theorem
\ref{critical-thm} regarding the analytical properties of the
critical time $t_c(K)$.

\subsection{Asymptotical Behavior of $t_c(K)$ for $K \gg 1$}
\begin{proof}[\textbf{\em Proof of Theorem \ref{critical-thm}}]
The proof begins with the analysis of $t_c(K)$ as $K\to\infty$,
whose leading order asymptotics is constructed by the method of
matched asymptotics.

The differential equation (\ref{y-eq}) for $y(x)$
\begin{equation}
y' = \frac{(1-y)K}{y^2+(1-y^2)K}-1,
\end{equation}
with its initial condition $y(0)=1$, suggests the transformation
$u=1-y$, so $u \ll 1$ near the origin, and satisfies
\[
-u' = \frac{Ku}{(1-u)^2 + Ku(2-u)}-1,
\]
with the initial condition $u(0)=0$. Hereafter $\approx$ means
equality of leading order terms. For $u\ll 1$ we have
\begin{equation}
\label{eq:u'} u' \approx \frac{1+Ku}{1+2Ku},
\end{equation}
whose solution is
\begin{equation}
\label{eq:u} 2Ku - \ln(1+Ku) \approx Kx.
\end{equation}
For $K^{-1} \ll u \ll 1$ the logarithmic term is smaller compared to
the linear term, hence
\begin{equation}
\label{eq:u-approx} u(x) \approx \frac{x}{2},\quad \mbox{for }
K^{-1} \ll x \ll K^{-\varepsilon},
\end{equation}
where $\varepsilon >0$ is needed to ensure $u\ll 1$.

We rewrite the differential equation \eqref{z-eq} for $z(x)$ as
\begin{equation}
\label{eq:riccati} z' = c(x)(z^2-1) + 1,
\end{equation}
where
\begin{equation}
\label{eq:c} c(x) = \frac{K}{y^2 + (1-y^2)K}.
\end{equation}
We recognize \eqref{eq:riccati} as a Riccati equation (see, e.g.,
\cite{Murphy}). However, it has no closed-form solution for a
general function $c(x)$, and we resort to asymptotic methods.

We obtain an asymptotic approximation of $z(x)$ by the method of
matched asymptotics. As shown below, we first construct an
asymptotic solution that is valid for $K^{-1} \ll x \ll K^{-2/3}$,
which is then matched to a different asymptotic solution that is
valid for $x \gg K^{-3/4}$ up to the blowup point.

First, we use a linear approximation of (\ref{eq:riccati}) near
$x=0$, by substituting $v=z-1$ that satisfies
\begin{equation}
v' = c(v^2+2v)+1,
\end{equation}
with the initial condition $v(0)=0$, so $v$ is small near the
origin. Therefore, as long as the nonlinear term satisfies $cv^2 \ll
1$, it can be viewed as a perturbation to the linear equation
\begin{equation}
\label{eq:diff-v} v' - 2cv = 1 + O(cv^2).
\end{equation}
Multiplying (\ref{eq:diff-v}) by its integrating factor we obtain
\begin{equation}
\label{eq:int-factor} v(x) \approx \exp\left\{\int 2c(x)\,dx
\right\} \int \exp\left\{-\int 2c(x')\,dx' \right\}\,dx,
\end{equation}
Equations (\ref{eq:u'}) and (\ref{eq:c}) imply
\begin{equation}
c \approx \frac{Ku'}{1+Ku}.
\end{equation}
Therefore,
\begin{equation}
\label{eq:int-c} \int c(x)\,dx \approx \int \frac{Ku'}{1+Ku}\,dx =
\ln (1+Ku).
\end{equation}
Hence, by equation (\ref{eq:int-factor}) we have
\begin{equation}
\label{eq:v-int} v \approx (1+Ku)^2\int \frac{dx}{(1+Ku)^2}.
\end{equation}
We calculate the integral in equation (\ref{eq:v-int}) by implying
equation (\ref{eq:u'})
\begin{eqnarray}
\int \frac{dx}{(1+Ku)^2} &=& \int
\frac{u'}{(1+Ku)^2}\frac{1+2Ku}{1+Ku}\,dx = \int
\left[\frac{2}{(1+Ku)^2}-\frac{1}{(1+Ku)^3} \right]du \nonumber \\
&=& \frac{1}{K}\left[-\frac{2}{1+Ku} + \frac{1}{2(1+Ku)^2} +
\frac{3}{2} \right] \nonumber
\end{eqnarray}
Thus,
\[
v \approx \frac{1}{2}u(2+3Ku),
\]
which combined with equation (\ref{eq:u-approx}) for $x \gg K^{-1}$
gives
\begin{equation}
\label{eq:v-approx} v(x) \approx \frac{3}{8}Kx^2.
\end{equation}
Equations (\ref{eq:c}) and (\ref{eq:u-approx}) imply that $c(x)$ is
to leading order
\begin{equation}
\label{eq:c1} c(x) \approx \frac{K}{1+2Ku} \approx \frac{1}{x},
\quad \mbox{for } K^{-1} \ll x \ll K^{-\varepsilon}.
\end{equation}
Recall that the linear approximation for $v(x)$ is valid whenever
$cv^2 \ll 1$. Equations (\ref{eq:v-approx}) and (\ref{eq:c1}) give
$cv^2 \approx K^2 x^3$, hence the linear approximation
(\ref{eq:diff-v}) holds for $x \ll K^{-2/3}$. Thus, the
approximation (\ref{eq:v-approx}) holds for $K^{-1} \ll x \ll
K^{-2/3}$.

Second, we choose an intermediate point $x_1$ such that $K^{-1} \ll
x_1 \ll K^{-2/3}$, e.g., $x_1 = K^{-3/4}$, for which
$c(x_1)=O\left(K^{3/4}\right)$, $v(x_1)=O\left(K^{-1/2} \right)$,
and $z(x_1)^2-1=O\left(K^{-1/2} \right)$. Therefore,
$c(x_1)\left(z(x_1)^2-1\right) = O\left(K^{1/4} \right)$ and we may
approximate equation (\ref{eq:riccati}) by neglecting the constant
term $1 \ll K^{1/4}$
\begin{equation}
z' \approx c(x)(z^2-1),
\end{equation}
with the initial condition
\begin{equation}
\label{initial} z_1\equiv z(x_1) = 1 + v(x_1) \approx
1+\frac{3}{8}Kx_1^2.
\end{equation}
The blowup point $t_c$ satisfies
\begin{equation}
\label{eq:int-eq} \int_{z_1}^\infty \frac{dz}{z^2-1} =
\int_{x_1}^{t_c} c(x)\,dx.
\end{equation}
The left hand side of (\ref{eq:int-eq}) is
$$\int_{z(x_1)}^\infty \frac{dz}{z^2-1} =
\frac{1}{2}\ln\left(\frac{z_1+1}{z_1-1} \right),$$ whereas the right
hand side of (\ref{eq:int-eq}) is calculated by equation
(\ref{eq:int-c}) and the approximation (\ref{eq:u-approx})
$$\int_{x_1}^{t_c} c(x)\,dx \approx \ln \left(\frac{1+Ku_c}{1+Ku_1} \right) \approx \ln \left(\frac{1+\frac{1}{2}Kt_c}{1+\frac{1}{2}Kx_1} \right).$$
Therefore,
\begin{equation}
\frac{z_1+1}{z_1-1} \approx
\left(\frac{1+\frac{1}{2}Kt_c}{1+\frac{1}{2}Kx_1} \right)^2 \approx
\frac{t_c^2}{x_1^2}
\end{equation}
Substituting (\ref{initial}) and rearranging we obtain
\begin{equation}
t_c^2 \approx x_1^2 \frac{2+\frac{3}{8}Kx_1^2}{\frac{3}{8}Kx_1^2}
\approx \frac{16}{3K}.
\end{equation}
Note that the result is independent of the choice of the
intermediate point $x_1$, as expected. Hence
\begin{equation}
t_c \approx \frac{4}{\sqrt{3K}}, \quad \mbox{for } K \gg 1,
\end{equation}
is the leading order asymptotics of the blow up point. Throughout
the asymptotic analysis we used merely the leading order terms in
their validity regimes. Therefore, the higher order terms are
guaranteed to be asymptotically smaller than the leading order term,
that is,
\begin{equation}
t_c(K) = \frac{4}{\sqrt{3K}}\left(1 + o(1)\right), \quad \mbox{for }
K \gg 1.
\end{equation}

We constructed numerical solutions of the ODEs
\eqref{y-eq},\eqref{z-eq} for various values of $K$ using
Mathematica. Figure \ref{fig::sing-approx} shows the comparison of
the numerical approximation for $t_c(K)$ and the asymptotic result
$4/\sqrt{3K}$.

\begin{figure}
\centering
\includegraphics{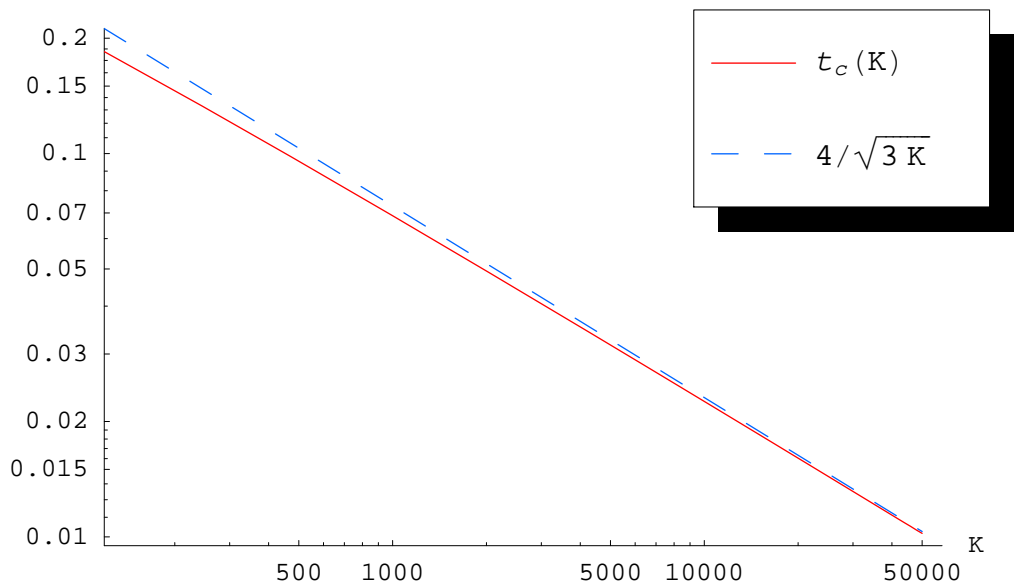}
\caption{Comparison of the numerical results for $t_c(K)$ obtained
by numerical solutions of the ODEs (\ref{y-eq})-(\ref{z-eq}), and
the asymptotic approximation of Theorem \ref{critical-thm}.
Logarithmic scale was used in both axes.} \label{fig::sing-approx}
\end{figure}

\subsection{Strict monotone decreasing behavior of $t_c(K)$}
The continuity of the critical point $t_c(K)$ follows from the
general continuous dependence of ODEs on their parameters, which is
in our case a single parameter $K$. Thus, to complete the proof we
must show that $t_c(K)$ is strictly monotone decreasing in $K$.

We first show that $c(x)$ is strict monotone increasing in $K$,
i.e., for any $0<K_1<K_2$, the corresponding functions $c_1,c_2$
(where $y$ is replaced by the corresponding solutions to
\eqref{y-eq}) satisfy $c_1(x)<c_2(x)$ for every $x \geq 0$.

Recall that by the properties of $y$, as stated in Proposition
\ref{prop-y-sol-exp-bound}, $c(0)=K$ and for every $K\neq 0,1$,
$c(x)$ either strictly increases or strictly decreases to $1$ as $x
\to \infty$ (depending on whether $K>1$ or $K<1$). Hence, if $K_1
\leq 1 < K_2$ or $K_1 < 1 \leq K_2$, clearly $c_1(x)<c_2(x)$ for
every $x \geq 0$. Furthermore, if $1 < K_1 < K_2$, the strict
monotonicity of $y$ in $K$ implies that for every $x > 0$:
$$c_1 = \frac{1}{y_1^2(\frac{1-K_1}{K_1})+1} <
\frac{1}{y_2^2(\frac{1-K_2}{K_2})+1} = c_2~.$$ Assume therefore that
$0 < K_1 < K_2 < 1$. Chain rule differentiation of \eqref{eq:c},
together with \eqref{y-eq}, gives:
\begin{equation}
\label{eq:c'-of-c-K-y}c'(x)=-2c^2\frac{1-K}{K}y\left((1-y)c+1\right)~.
\end{equation}
Rearranging equation \eqref{eq:c}, we get:
\begin{equation}
\label{eq:y-of-c} y = \left(\frac{(1-c)K}{c(1-K)}\right)^{1/2}.
\end{equation}
Combining equations \eqref{eq:c'-of-c-K-y} and \eqref{eq:y-of-c}, we
obtain an equality of the form $c'(x) = g(c,K)$, where $g(c,K)$
satisfies:
\begin{equation}
\label{eq:d-c'-over-d-K}\frac{\partial g}{\partial K}(c,K) =
\frac{c^2 (1+c)}{K} \sqrt{\frac{(1-c)K}{c(1-K)}} > 0~,
\end{equation}
for all $0 < K < 1$ and $x \geq 0$. By the argument stated in the
proof of Proposition \ref{prop-y-sol-exp-bound}, this implies that
$c$ is indeed strictly monotone increasing in $K$.

Returning to equation \eqref{eq:riccati}, since $z' = c(z^2-1)+1$
and $c$ is strictly monotone increasing in $K$, we conclude that
$z'(t)$ is monotone increasing in $K$ for any $t>0$ (recall that, as
stated in the proof of Proposition \ref{z-thm}, $z'(t)>0$ for any $t
\geq 0$ and hence $z(t)> 1$ for any $t>0$). Hence, the above
consideration implies that $z$ is strictly monotone increasing in
$K$. Therefore, $t_c(K)$ is monotone decreasing in $K$, that is,
\begin{equation}
\label{eq:mono} t_c(K_1) \geq t_c(K_2), \quad \mbox{for } K_1 < K_2.
\end{equation} Yet, the fact that $z$ is strictly monotone increasing in $K$ does
not guarantee by itself that $t_c(K)$ is {\em strictly} monotone
decreasing in $K$.

To this end, we consider solutions $z_1,z_2$ and $c_1,c_2$
corresponding to $K_1 < K_2$ and prove that $t_c(K_1) > t_c(K_2)$.
The functions $c_1$ and $c_2$ are continuous on $[0,\infty)$ and
satisfy $c_1 < c_2$, therefore $\exists \varepsilon > 0$ such that
\begin{equation}
\label{eq:c1c2} c_2(x) - c_1(x) > \varepsilon, \;  \forall x\in
[0,t_c(K_1)].
\end{equation}
The function (whose choice is made clear below) $$A(z) =
\frac{\ln\ds\frac{z^2}{z^2-1}}{\ln \ds\frac{z+1}{z}}$$ satisfies
$\lim_{z\to\infty} A(z) = 0$. Since $z_1(x)\to \infty$ as $x \to
t_c(K_1)$, there exists a point $0 <t^* < t_c(K_1)$, such that
\begin{equation}
\label{eq:conds} z_1(t^*) \geq \max \{1, K_1^{-1} \} \quad
\mbox{and} \quad A(z_1(t^*)) < \frac{2\varepsilon}{\max\{1,K_1\}}.
\end{equation}
Let $z_1^*(x)$ be the solution of the differential equation
\begin{equation}
\label{eq:z1*} z' = c_1(x)z^2 + c_1(x)z, \quad z(t^*) = z_1(t^*).
\end{equation}
Recall that $z_1(x)$ is the solution of
\begin{equation}
\label{eq:z1} z' = c_1(x)z^2 + 1-c_1(x), \quad z(t^*) = z_1(t^*).
\end{equation}
Comparing (\ref{eq:z1*}) and (\ref{eq:z1}) we note that
$c_1(x)z_1^*(x)
> 1-c_1(x)$ for all $x\geq t^*$, because $z_1^*(x)$ is monotone
increasing in $x$, $c_1(x) \geq \min\{1,K_1 \}$, and the choice of
$t^*$ in equation (\ref{eq:conds}), indicated by
$$c_1(x)z_1^*(x) \geq c_1(x)z_1^*(t^*) \geq
\min\{1,K_1\} \max\{1,K_1^{-1}\} = 1 > 1-c_1(x).$$ Hence, $z_1^*(x)
> z_1(x)$ for all $x
> t^*$ and
\begin{equation}
\label{eq:t1*} t_c(K_1) \geq t_c^*(K_1),
\end{equation}
where $t_c^*(K_1)$ is the blowup point of $z_1^*(x)$. Similar
considerations show that the blowup point $t_c^*(K_2)$ of the
solution, $z_2^*$, to the differential equation
\begin{equation}
\label{eq:z2*} z' = c_2(x)(z^2-1), \quad z(t^*) = z_2(t^*),
\end{equation}
(remark: if $t^* \geq t_c(K_2)$ then $t_c(K_1) > t_c(K_2)$ and the
proof is completed, so we assume $t^* < t_c(K_2)$) satisfies
\begin{equation}
\label{eq:t2*} t_c^*(K_2) \geq t_c(K_2),
\end{equation}
because the positive term 1 was omitted in (\ref{eq:z2*}) compared
to (\ref{eq:riccati}). The advantage of the differential equations
(\ref{eq:z1*}) and (\ref{eq:z2*}) over the original equation
(\ref{eq:riccati}) is that their blowup points satisfy closed form
relations. Namely,
\begin{eqnarray}
\label{eq:z1t1} \ln \left(\frac{z_1^*(t^*)+1}{z_1^*(t^*)} \right)
&=&
\int_{t^*}^{t_c^*(K_1)}c_1(x)\,dx, \\
\label{eq:z2t2} \frac{1}{2}\ln
\left(\frac{z_2^*(t^*)+1}{z_2^*(t^*)-1} \right) &=&
\int_{t^*}^{t_c^*(K_2)}c_2(x)\,dx.
\end{eqnarray}
Assume to the contrary that $t_c(K_1)=t_c(K_2)$. Then, by equations
(\ref{eq:t1*}) and (\ref{eq:t2*}) it follows that
\begin{equation}
\label{eq:t1t2} t_c^*(K_2) \geq t_c^*(K_1).
\end{equation}
Combining equations
(\ref{eq:z1t1}),(\ref{eq:z2t2}),(\ref{eq:t1t2}),(\ref{eq:c1c2}) and
the positivity of $c_2(x)$, we obtain
\begin{eqnarray}
\frac{1}{2}\ln \left(\frac{z_2^*(t^*)+1}{z_2^*(t^*)-1} \right)
&\geq& \int_{t^*}^{t_c^*(K_1)}c_2(x)\,dx > \varepsilon
(t_c^*(K_1)-t^*) + \int_{t^*}^{t_c^*(K_1)}c_1(x)\,dx \nonumber \\
&&\\
&=& \ln \left(\frac{z_1^*(t^*)+1}{z_1^*(t^*)} \right) + \varepsilon
(t_c^*(K_1)-t^*). \nonumber
\end{eqnarray}
Since $\ln\frac{z+1}{z-1}$ is a monotonic decreasing function of
$z$, and $z_1^*(t^*) < z_2^*(t^*)$, it follows that
\begin{equation}
\label{eq:bound1} \ln \left(\frac{z_1^*(t^*)+1}{z_1^*(t^*)-1}
\right) - 2\ln \left(\frac{z_1^*(t^*)+1}{z_1^*(t^*)} \right) >
2\varepsilon (t_c^*(K_1)-t^*),
\end{equation}
or equivalently
\begin{equation}
\ln \left(\frac{z_1^*(t^*)^2}{z_1^*(t^*)^2-1}\right) > 2\varepsilon
(t_c^*(K_1)-t^*).
\end{equation}
Employing equation (\ref{eq:z1t1}) with $c_1(x) \leq \max\{1,K_1\}$,
we find
\begin{equation}
\label{eq:bound2} t_c^*(K_1)-t^* \geq \frac{1}{\max\{1,K_1\}} \ln
\left(\frac{z_1^*(t^*)+1}{z_1^*(t^*)} \right).
\end{equation}
Equations (\ref{eq:bound1}) and (\ref{eq:bound2}) give
\begin{equation}
A(z_1^*(t^*)) = \frac{\ln
\left(\ds\frac{z_1^*(t^*)^2}{z_1^*(t^*)^2-1}\right)}{\ln
\left(\ds\frac{z_1^*(t^*)+1}{z_1^*(t^*)} \right)} >
\frac{2\varepsilon}{\max\{1,K_1\}},
\end{equation}
which contradicts (\ref{eq:conds}). Therefore, $t_c(K_1) >
t_c(K_2)$.
\end{proof}

\section{Concluding remarks and open
problems}\label{sec::conclusion} As mentioned in
\cite{SpencerWormald}, there is a well known relation between the
``double jump'' phenomenon in $\Gorg[1]$ and the critical threshold
in percolation. Indeed, the susceptibility corresponds to the size
of the component containing the origin, and the appearing of giant
component is analogous to the appearing of an infinite cluster. The
subcritical and supercritical stages in the $\Gorg[1]$ correspond to
those in percolation. However, it is unclear how to define a
percolation process which is analogous to the biased process.

We conclude by describing several natural ways to extend the family
of biased processes and their analysis.

Consider processes where the probability of adding an edge at each
step is biased according to the sizes of the components of its
endpoints. For instance, it is possible to assign the weight
$K_{i,j}$ to a missing edge $(u,v)$, such that
$|\mathcal{C}_u|=i$,$|\mathcal{C}_u|=j$ for $i,j<C$, where $C>0$ is
some large constant, and some weight $K_C$ otherwise. It seems
plausible that the above methods can be applied in this scenario
as-well.

The biased process $\Gorg$ was defined as a graph process on the
complete graph on $n$ vertices, $K_n$; in other words, at each step
we considered the entire set of missing edges, and the process ended
when obtaining $K_n$. It is possible to define an analogous process
on a different underlying graph $H=(V,E)$, where the edge pool at
each step contains only edges belonging to $E$, and the process ends
once $H$ is obtained. For instance, it is possible to consider the
biased process on $Q_n$, the $n$-th dimensional cube, or $L_2(n)$,
the $n\times n$ grid. In these cases, it may be interesting to study
the critical time for satisfying other monotone increasing
properties, such as left-right crossing in $L_2(n)$. It seems
plausible that this crossing probability is monotone decreasing in
$K$.

\noindent\textbf{Acknowledgement} We would like to thank Itai
Benjamini, who proposed the model of the biased process and provided
useful insights to the problem.

\end{document}